\documentclass[pdflatex,sn-mathphys-num]{sn-jnl}

\usepackage{graphicx}%
\usepackage{multirow}%
\usepackage{amsmath,amssymb,amsfonts}%
\usepackage{amsthm}%
\usepackage{mathrsfs}%
\usepackage[title]{appendix}%

\usepackage{textcomp}%
\usepackage{manyfoot}%
\usepackage{booktabs}%
\usepackage{algorithm}%
\usepackage{algorithmicx}%
\usepackage{algpseudocode}%
\usepackage{listings}%

\usepackage{hyperref}
\usepackage{xcolor}

\hypersetup{
	colorlinks=true,
	urlcolor=blue, 
	citecolor=blue 
}
\usepackage{subcaption}
\usepackage{float}

\theoremstyle{thmstyleone}%

\newtheorem{theorem}{Theorem}[section] 
\newtheorem{lemma}{Lemma}[section]

\begin{document}
	
	\title[Article Title]{3-path-connectivity of bubble-sort star graphs}
	
	
	\author[1]{Yi-Lu Luo}
	
	\author*[1]{Yun-Ping Deng}\email{dyp612@hotmail.com}

	\author[1]{Yuan Sun}

	\affil[1]{\orgdiv{Department of Mathematics}, \orgname{Shanghai University of Electric Power}, \orgaddress{\street{2103 Pingliang Road}, \city{Shanghai}, \postcode{200090}, \country{PR China}}}

	
	\abstract{	Let $G$ be a simple connected graph with vertex set $V(G)$ and edge set $E(G)$. Let $T$ be a subset of $ V(G)$ with cardinality $|T|\geq2$. A path connecting all vertices of $T$ is called a $T$-path of $G$. Two $T$-paths $P_i$ and $P_j$ are said to be internally disjoint if $V(P_i)\cap V(P_j)=T$ and $E(P_i)\cap E(P_j)=\emptyset$. Denote by $\pi_G(T)$ the maximum number of internally disjoint $T$-paths in G. Then for an integer $\ell$ with $\ell\geq2$, the $\ell$-path-connectivity  $\pi_\ell(G)$ of $G$ is formulated as $\min\{\pi_G(T)\,|\,T\subseteq V(G)$ and $|T|=\ell\}$. In this paper, we study the $3$-path-connectivity of  $n$-dimensional bubble-sort star graph $BS_n$. By deeply analyzing the structure of $BS_n$, we show that  $\pi_3(BS_n)=\lfloor\frac{3n}2\rfloor-3$, for any $n\geq3$.}

	
	
	
	\maketitle
	\section{Introduction}\label{sec1}
	With the rapid advancement of technology, the design and analysis of various interconnection networks have become an important area of research. The underlying topology of interconnection networks is usually modeled as a simple connected graph. Connectivity is a crucial parameter for evaluating the reliability and fault tolerability of a graph.
	Let $G$ be a simple connected graph with vertex set $V(G)$ and edge set $E(G)$. 
	For any $2$-subset $\{u, v\} \subseteq V(G)$, let $\kappa_G(u, v)$ denote the maximum number of
	internally disjoint paths connecting vertices $u$ and $v$ in $G$.
	The veterx connectivity $\kappa(G)$ of $G$ is defined as $\min\{\kappa_G(u, v)\mid\{u, v\} \subseteq V(G)\}$.
	This is a specialized definition of vertex connectivity, proposed by Whitney \cite{whitney1992congruent}, which emphasizes the maximum number of internally disjoint paths between any two vertices in a graph.
	As a generalization of vertex connectivity, Hager \cite{Hager86} proposed the concept of $\ell$-path-connectivity, which measures the fault tolerance of a graph by the maximum number of internally disjoint paths connecting any $\ell$ vertices in the graph.
	Let $T$ be a subset of $V(G)$ with cardinality $|T|\geq2$. A path connecting all vertices of $T$ is called a $T$-path of $G$.
	Two $T$-paths $P_i$ and $P_j$ are said to be internally disjoint if $V(P_i)\cap V(P_j)=T$ and $E(P_i)\cap E(P_j)=\emptyset$.
	Denote by $\pi_G(T)$ the maximum number of internally disjoint $T$-paths in $G$. Then for an integer $\ell$ with $\ell\geq2$, the $\ell$-path-connectivity  $\pi_\ell(G)$ of $G$ is formulated as $\min\{\pi_G(T)|T\subseteq V(G)$ and $|T|=\ell\}$.
	In fact, $\pi_2(G)$ is the vertex connectivity $\kappa(G)$.

	There are some known results about path-connectivity. For example, Mao \cite{Mao16} determined the upper and low bounds of $3$-path-connectivity for the lexicographic product of two connected graphs.
	Li et al. \cite{Li21} proved that deciding whether $\pi_G(T)\geq k$ with $T \subseteq V(G)$ is NP-complete for any $k \geq1$.
	The $3$-path-connectivity of some graphs has been completely determined, such as hypercubes \cite{Zhu22}, $k$-ary $n$-cube \cite{Zhu23}, Cayley graphs generated by transposition trees \cite{Jin23}, star graphs \cite{Li24} and pancake graphs \cite{wang25}.
	
	Cayley graphs, especially on the symmetric groups, are often
	used as the models of interconnection networks \cite{Akers89,Lakshmivarahan93}.
	Let $X$ be a finite group and $S$ an inverse-closed subset of $X \setminus \{e\}$, where $e$ is the identity element of $X$.
	The Cayley graph on $X$ with respect to $S$ is the graph with vertex set $X$ and edge set $\{(g, gs) \,|\, g\in X, s\in S\}$. Let  $S_n$ be the symmetric group of degree $n$ with $n\geq3$.
	For convenience, we denote by $(x_1x_2\ldots x_n)$ the permutation 
	$\left(\begin{array}
		{llll}1 & 2 & \ldots & n \\
		x_1 & x_2 & \ldots & x_n
	\end{array}\right)$ and by $(i,j)$ the transposition that exchanges the digits at positions $i$ and $j$, that is, $(x_1x_2\ldots x_i\ldots x_j\ldots x_n)\circ(i,j)=(x_1x_2\ldots x_j\ldots x_i\ldots x_n)$. Set $S=\{(1,j)\mid2\leqslant j\leqslant n\}\cup\{(j,j+1)\mid2\leqslant j\leqslant n-1\}$. 
	The bubble-sort star graph, denoted by $BS_n$, is defined as the Cayley
	graph on $S_n$ with respect to $S$. 
	This graph is introduced as an interconnection network topology in \cite{Chou96},  and several variants of its connectivity have been studied \cite{Guo16,Guo21,Gu20,Zhu18,Zhang19,cheng2023path,zhao19}. 
	
	Motivated by the above results, in this paper, we focus on the 3-path-connectivity of $BS_n$. In Section \ref{sec2}, we obtain some properties of $BS_n$ and lemmas that will be used in the following sections. In Section \ref{sec3}, we obtain a structure connecting any three vertices in $V(BS_n)$. The main result $\pi_3(BS_n)=\lfloor\frac{3n}2\rfloor-3$ with  any $n\geq3$ is derived in Section \ref{sec4}.

	\section{Preliminaries}
	\label{sec2}
	Let $G$ be a simple connected graph and $H$ a subgraph of $G$. Set $G\setminus H=G[V(G)-V(H)]$, where $G[V(G)-V(H)]$ is a subgraph of $G$ induced by $V(G)-V(H)$. For any $u\in V(G)$, the set of neighbors of $u$ in $G$ is denoted by $N_G(u)$. For $\{u_1,\ldots,u_k\}\subseteq V(G)$, the set of common neighbors of $u_1,\ldots,u_k$ in $G$ is denoted by $CN_G(u_1,\ldots,u_k)$, that is, $CN_G(u_1,\ldots,u_k)=N_{G}(u_1)\cap\ldots\cap N_{G}(u_k)$.
	We denote $\{1,2,\ldots,n\}$ by $[n]$. 
	Regarding the notations not defined in this paper, we
	refer the readers to \cite{Bondy08}.
	
	Let $BS_n^i$ be the subgraph of $BS_n$ induced by the set $\{(x_1x_2\ldots x_n)\in S_n\mid x_n=i\}$ for $i\in[n]$. Clearly, each $BS_n^i$ is isomorphic to $BS_{n-1}$, forming a substructure copy of $BS_n$.
	Thus, $BS_n$ can decompose into $n$ such copies, denoted as $BS_n = BS_n^1 \oplus BS_n^2 \oplus \cdots \oplus BS_n^n$, where $\oplus$ represents the decomposition relationship.
	It is easy to see that for any $u\in V(BS_n^i)$,  $u$ has exactly $2n-5$ neighbors within $BS_n^i$
	and two neighbors outside $BS_n^i$, which are called the outgoing neighbors of
	$u$ and denoted by $u^+ = u \circ ( 1, n)$, $ u^- = u \circ ( n- 1, n)$. Clearly, $u^+$ and $u^-$ belong to different copies, and  $\{u^+,u^-\}\cap\{v^+,v^-\}=\emptyset$ for any two distinct vertices $u, v\in V(BS_n^{i})$.

	\begin{lemma}\label{lem2.1} (\cite{Cai15})
		For $i,j\in [n]$ with $i\neq j$, the following properties hold:
		
		(i) $BS_{n}$ is $2n-3$-regular and $\kappa(BS_{n})=2n-3$ for $n\geq2$.
		
		(ii)  For $n\geq 4$, $\left|E_{i, j}(BS_{n})\right|=2(n-2)!$, where  $E_{i, j}(BS_n)$ is the set of edges between $BS_n^i$ and $BS_n^j$.
		
		(iii) For $n\geq3$ and $\{u, v\}\subseteq V(BS_n)$, $\left|CN_{BS_n}(u,v)\right|\leq 3$.

	\end{lemma}

	\begin{figure}[H]
		\centering
		\begin{subfigure}[b]{0.4\textwidth} 
			\includegraphics[width=\textwidth]{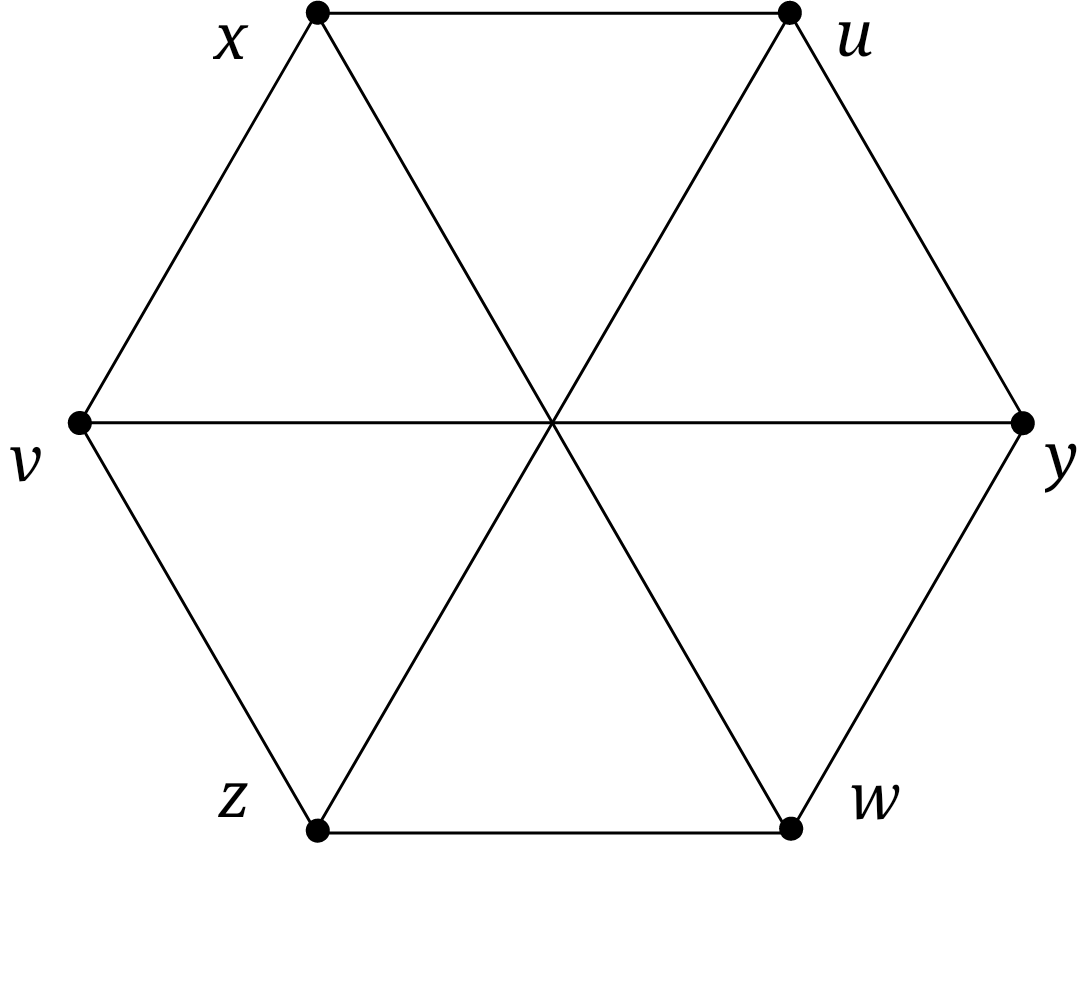} 
			\subcaption{$BS_3$}
			
		\end{subfigure}
		\hspace{35pt}
		\begin{subfigure}[b]{0.4\textwidth} 
			\includegraphics[width=\textwidth]{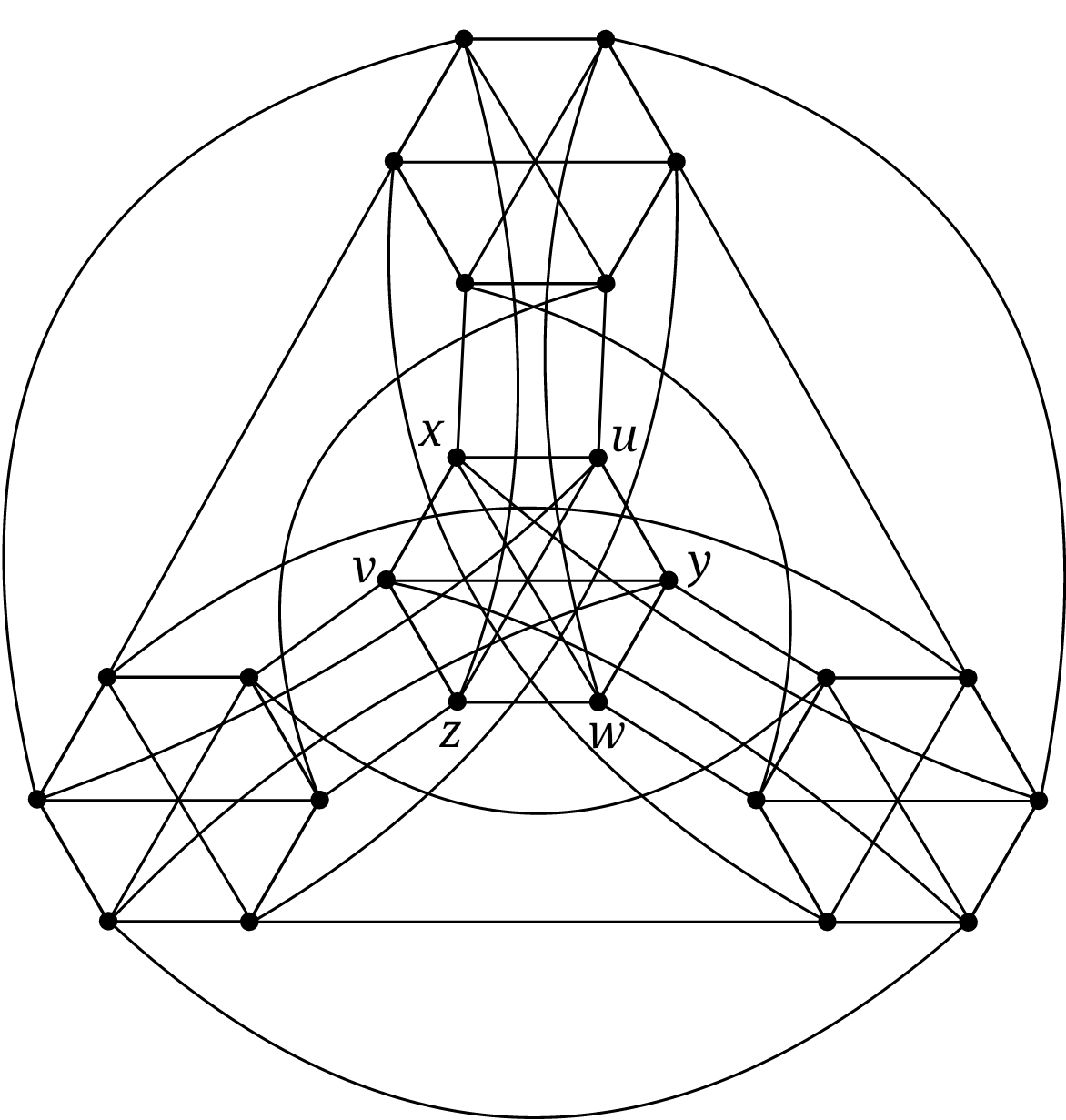} 
			\subcaption{$BS_4$}
			
		\end{subfigure}
		\caption{The bubble-sort star graphs $BS_3$ and $BS_4$ }
		\label{fig1}
	\end{figure}
	
	\begin{lemma}\label{lem2.2} 
		For any $n\geq 3$, $max\left \{ \left | CN_{BS_n}( u,v,w) \right | \mid \{ u,v,w\} \subseteq V(BS_n) \right \} =3$.
		\\
		\\
		\noindent\textbf{Proof.} 
		By Fig. \ref{fig1} (a), we have  $\left|CN_{BS_3}(u,v,w)\right|=\left |\{x,y,z\}\right |=3$. Since $BS_3$ is a subgraph of $BS_n$ with $n\geq3$, it follows that  $max\left \{ \left | CN_{BS_n}( u,v,w) \right | \mid \{ u,v,w\} \subseteq V(BS_n) \right \} \geq3$.
		
		On the other hand, by Lemma \ref{lem2.1} (iii), $\left|CN_{BS_n}(u,v)\right|\leq3$, and thus $\left|CN_{BS_n}(u,v,w)\right|$ $\leq3$ for any $u,v,w\in V(BS_n)$, that is, $max\left \{\left|CN_{BS_n}(u,v,w)\right|\mid\{u,v,w\}\subseteq V(BS_n)\right \} {\leq}3$.
		
		
		Therefore, $max\left \{ \left | CN_{BS_n}( u,v,w) \right | \mid \{ u,v,w\} \subseteq V(BS_n) \right \}=3$. $\square$
	\end{lemma}

	\begin{lemma}\label{lem2.3} (\cite{Bondy08})
		Let $G$ be a $k$-connected graph, let $x$ be a vertex of $G$, and let $Y \subseteq V(G) \setminus \{x\}$ be a set of at least $k$ vertices of $G$. Then there exists a $k$-fan in $G$ from $x$ to $Y$.
	\end{lemma}
	
	\begin{lemma}\label{lem2.4} (\cite{Bondy08})
		Let $G$ be a $k$-connected graph, and let $X$ and $Y$ be subsets of $V(G)$ of cardinality at least $k$. Then there exists in $G$ a family of $k$ pairwise disjoint $(X,Y)$-paths.
	\end{lemma}
	
	\begin{figure}[H]
		\centering
		\begin{subfigure}[b]{0.4\textwidth} 
			\includegraphics[width=\textwidth]{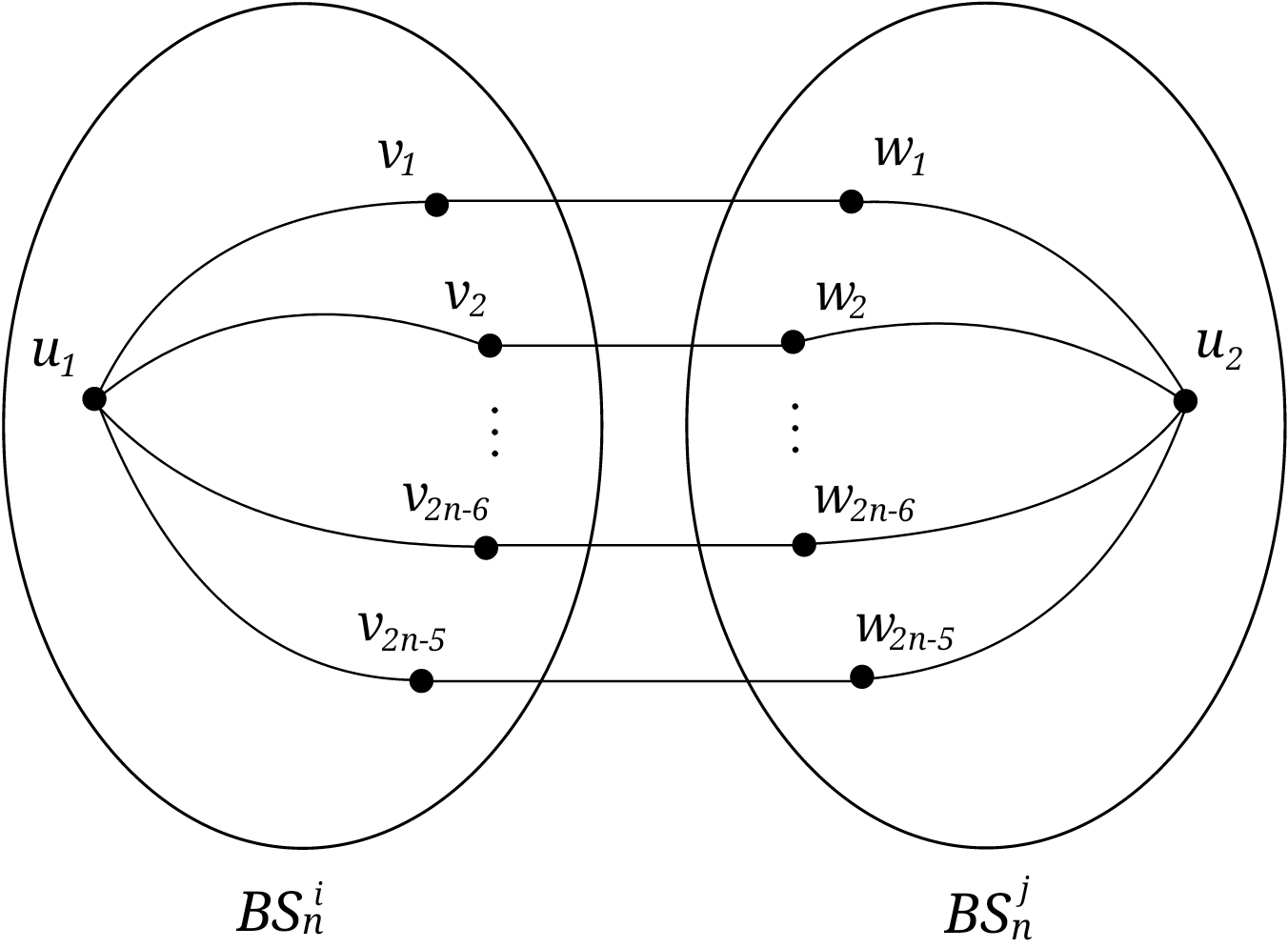} 
			\subcaption{Case 2 of (i)}
			
		\end{subfigure}
		\hspace{35pt}
		\begin{subfigure}[b]{0.4\textwidth} 
			\includegraphics[width=\textwidth]{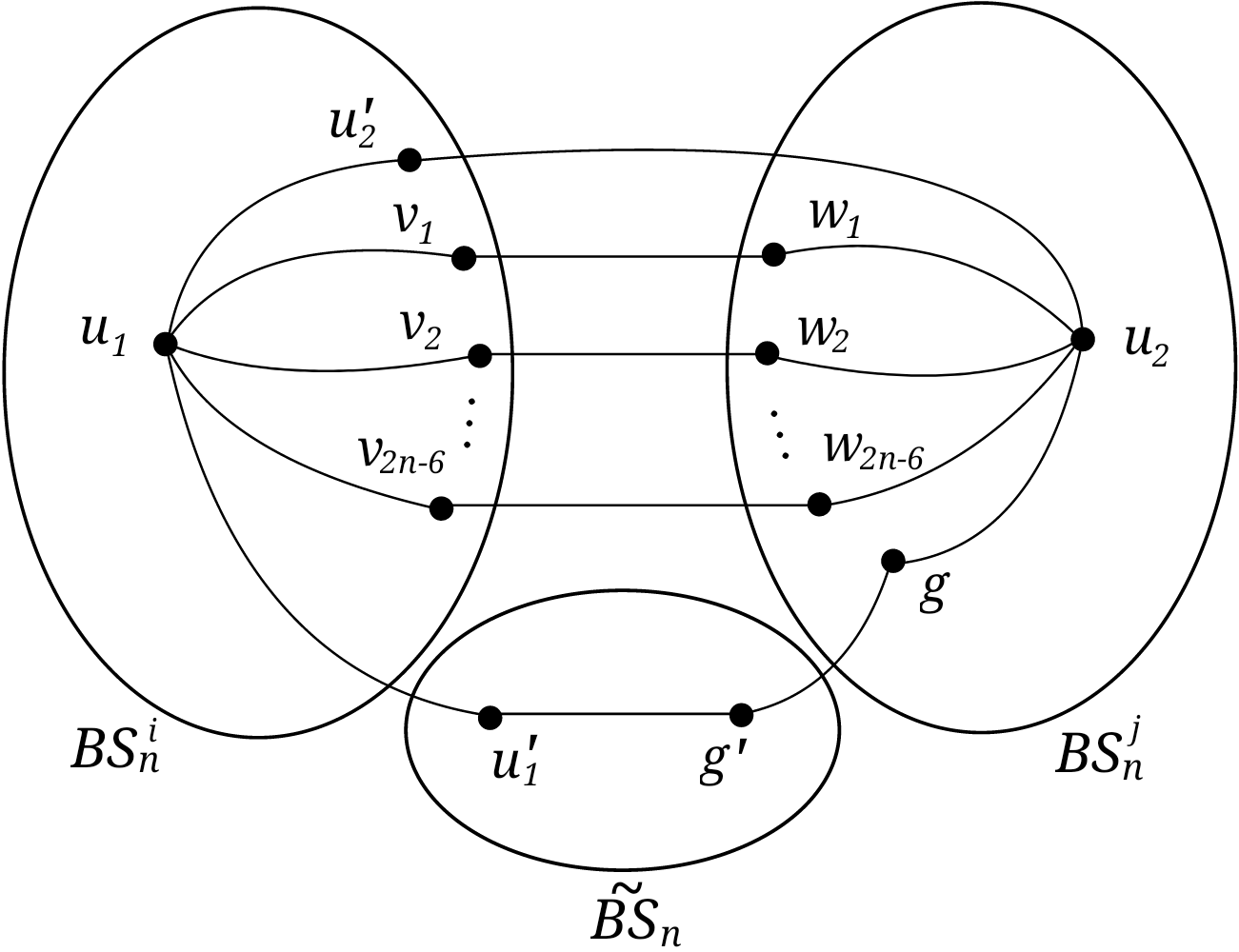} 
			\subcaption{Subcase 2.2 of (ii)}
			
		\end{subfigure}
		\caption{The illustrations of Lemma 2.5}
		\label{fig2}
	\end{figure}
	
	\begin{lemma}\label{lem2.5}
		For $n\geq3$, let $\widehat{BS_n}=BS_{n}^{t_{1}}\oplus BS_{n}^{t_{2}}\oplus\cdots\oplus BS_{n}^{t_{k}}$ for any $\{t_1,\ldots,t_k\}\subseteq\{1,2,\cdots,n\}$ with $1\leq k\leq n-1$, where $BS_{n}^{t_{1}}\oplus BS_{n}^{t_{2}}\oplus\cdots\oplus BS_{n}^{t_{k}}$ $is$ the induced subgraph of $BS_n$ on $V\left (BS_{n}^{t_{1}}\right ) \cup V\left (BS_{n}^{t_{2}}\right ) \cup \cdots \cup V\left ( BS_{n}^{t_k}\right )$. Then the  following properties hold:

		(i) For $1\leq k\leq n-2$, $\kappa(\widehat{BS_n})=2n-5$.

		(ii) For $k=n-1$, $\kappa(\widehat{BS_n})=2n-4$.
		\\
		\\
		\noindent\textbf{Proof.}
		\textbf{(i)} For $k=1$, by Lemma \ref{lem2.1} (i), $\kappa(\widehat{BS_n})=\kappa({BS_{n}^{t_1}})=2(n-1)-3=2n-5$.
		
		For $2\leq k\leq n-2$, there must exist a vertex $u \in V(\widehat{BS_n})$ such that $\{u^+,u^-\}$ is outside of $\widehat {BS_n}$, hence  $\delta({\widehat{BS_n}})=2n-5$, and it follows that $\kappa(\widehat{BS_n})\leq \delta(\widehat{BS_n})=2n-5$.

		Next, we need to prove that $\kappa(\widehat{BS_n})\geq 2n-5$. 
		
		For $n=3$, since $\widehat {BS_n}$ is connected, $\kappa(\widehat {BS_n})\geq1=2n-5$.
		
		For $n\geq4$, let $\{u_1 ,u_2\}$ be an arbitrary subset of $V(\widehat {BS_n})$. We prove that $\kappa(\{u_1,u_2\})\geq 2n-5$, where $\kappa(\{u_1,u_2\})$ is the maximum number of
		internally disjoint $(u_1 ,u_2)$-paths in $G$, and hence $\kappa(\widehat {BS_n})\geq 2n-5$. Then, we deliberate on the two cases that follow:

		\textbf{Case 1:} $u_1,u_2$ belong to the same copy.
		
		Without loss of generality, let $u_1,u_2\in V(BS_{n}^i)$. Since $\kappa(BS_{n}^i)=2n-5$, we can get $\kappa(\{u_{1},u_{2}\})\geq 2n-5$.
		
		\textbf{Case 2:} $u_1,u_{2}$ belong to two different copies.
		
		Assume that $u_1\in V(BS_{n}^i)$ and $u_2\in V(BS_{n}^j)$. Since $2n-5\leq2(n-2)!=\left|E_{i,j}(BS_n)\right| $, referring to Fig \ref{fig2} (a), we can select $2n-5$ edges $(v_t,w_t)$ in $E_{i,j}(BS_n)$, where $v_t\in V(BS_{n}^i)$ and $w_t\in V(BS_{n}^j)$ with $1\leq t\leq 2n-5$.
		
		Let $S_1=\{v_1,v_2,\ldots,v_{2n-5}\}$ and $S_2=\{w_1,w_2,\ldots,w_{2n-5}\}$.
		Recall that $\kappa(BS_{n}^i)$ $=2n-5$. By Lemma \ref{lem2.3}, there are $2n-5$ internally disjoint $(u_1,S_1)$-paths $P_1, P_2,\ldots, P_{2n-5}$ in $BS_{n}^i$ and the terminal vertex of $P_{i}$ is $v_{i}$ $(i\in[2n-5])$.
		Similarly, there are $2n-5$ internally disjoint $(u_2, S_2)$-paths $R_1, R_2,\ldots, R_{2n-5}$ in $BS_{n}^j$ and the terminal vertex of $R_{i}$ is $w_{i}$ $(i\in[2n-5])$.
		
		Now, we have obtained $2n-5$ internally disjoint $(u_1,u_2)$-paths: $u_1 P_1 v_1 w_1 R_1 u_2$, $u_1 P_2 v_2 w_2 R_2 u_2$,  $\ldots$ , $u_1 P_{2n-5} v_{2n-5} w_{2n-5} R_{2n-5} u_2$, that is, $\kappa(\{u_1,u_2\})\geq 2n-5$.
		
		In conclusion, $\kappa(\widehat {BS_n})=2n-5$.
		\\
		
		\textbf{(ii)} For $k=n-1$, for any vertex $u \in V(\widehat{BS_n})$, there is at most one outgoing neighbor $u^{\prime}$ is outside of $\widehat {BS_n}$, hence $\delta(\widehat{BS_n})=2n-4$, and it follows that $\kappa(\widehat{BS_n})\leq \delta(\widehat{BS_n})=2n-4$.
		
		Next, we need to prove that $\kappa(\widehat{BS_n})\geq 2n-4$. Let ${\widehat{BS_n}}=BS_n\setminus BS_{n}^s$ with $1\leq s\leq n$.
		
		For $n=3$, by Fig \ref{fig1} (a), $\kappa(\widehat {BS_n})\geq2=2n-4$.
		
		For $n\geq4$, let $\{u_1 ,u_2\}$ be an arbitrary subset of $V(\widehat {BS_n})$.
		By considering the following two cases, we prove that $\kappa(\{u_1,u_2\})\geq 2n-4$, and hence $\kappa(\widehat {BS_n})\geq 2n-4$.
		
		\textbf{Case 1:} $u_1,u_2$ belong to the same copy.
		
		Without loss of generality, let $u_1,u_2\in V(BS_{n}^i)$. Since $\kappa(BS_{n}^i)=2n-5$, there are $2n-5$ internally disjoint $(u_1,u_2)$-paths in $BS_{n}^i$.
		In addition, since $u_1$ and $u_2$ have at least one outgoing neighbor $u_1^{\prime}\notin V(BS_{n}^s)$ and $u_2^{\prime}\notin V(BS_{n}^s)$, respectively. Since   $\kappa(\widehat {BS_n}\setminus BS_{n}^i)\geq 1$, we can obtain another $(u_1,u_2)$-path $u_1 u_1^{\prime}u_2^{\prime}u_2$, whose internal vertices are in $\widehat {BS_n}\setminus BS_{n}^i$.
		Thus, $\kappa(\{u_1,u_2\})\geq 2n-4$.
		
		\textbf{Case 2:} $u_1,u_{2}$ belong to two different copies.
		
		Suppose that $u_1\in V(BS_{n}^i)$ and $u_2\in V(BS_{n}^j)$. 
		
		Let $\widetilde{BS_n}=\widehat {BS_n}\setminus({BS_{n}^i\oplus BS_{n}^j})=BS_n\setminus({BS_{n}^i\oplus BS_{n}^j}\oplus BS_{n}^s)$.
		
		\textbf{Subcase 2.1:} Both $u_1$ and $u_2$ have at least one outgoing neighbor in $\widetilde{BS_n}$.
		
		Let $u_1^{\prime}$ and $u_2^{\prime}$ be the outgoing neighbors of $u_1$ and $u_2$, respectively. Since $\kappa(\widetilde{BS_n})\geq 1$, we can obtain one $(u_1,u_2)$-path $u_1 u_1^{\prime}u_2^{\prime}u_2$, whose internal vertices are in $\widetilde{BS_n}$.
		
		Similarly to the proof of (i) Case 2, we can obtain another $2n-5$ internally disjoint $(u_1,u_2)$-paths in ${BS_{n}^i\oplus BS_{n}^j}$. Thus, $\kappa(\{u_1,u_2\})\geq 2n-4$.
		
		\textbf{Subcase 2.2:} Only $u_1$ has at least one outgoing neighbor in $\widetilde{BS_n}$.
		
		Let $u_1^{\prime}$ be the outgoing neighbor of $u_1$ in $\widetilde{BS_n}$ and $u_2^{\prime}$ be the outgoing neighbor of $u_2$ in $BS_{n}^i$. We can find a vertex $g\in V(BS_{n}^j)$ that has an outgoing neighbor $g^{\prime}\in V(\widetilde{BS_n})$. Since $\kappa(\widetilde{BS_n})\geq 1$, there exists a $(u_1^\prime,g^\prime)$-path in $\widetilde{BS_n}$.  
		
		Similarly to the proof of (i) Case 2, we can obtain a $2n{-}5$-fan from $u_1$ to $S_1\cup \{u_2^\prime\}$ in $BS_n^i$ and a $2n{-}5$-fan from $u_2$ to $S_2\cup \{g\}$ in $BS_n^j$.
		Therefore, referring to Fig \ref{fig2} (b) , we can obtain $2n-
		4$ internally disjoint $(u_1,u_2)$-paths: $u_1 P_1 v_1 w_1 R_1 u_2$, $u_1 P_2 v_2 w_2 R_2 u_2$, $\ldots$ , $u_1 P_{2n-6} v_{2n-6} w_{2n-6} R_{2n-6} u_2$, $u_1 u_1^{\prime} g^{\prime} g u_2,u_1 u_2^{\prime} u_2$. Thus, $\kappa(\{u_1,u_2\})\geq 2n-4$.

		\textbf{Subcase 2.3:} Neither  $u_1$ nor $u_2$ have any outgoing neighbor in $V(\widetilde{BS_n})$.
		
		Let $u_1^{\prime}$ be the outgoing neighbor of $u_1$ in ${BS_n^j}$ and $u_2^{\prime}$ be the outgoing neighbor of $u_2$ in $BS_{n}^i$.
		By Lemma \ref{lem2.3}, similar to the proof of (i) Case 2, we can obtain  $2n-4$ internally disjoint $(u_1,u_2)$-paths: $u_1 P_1 v_1 w_1 R_1 u_2$, $u_1 P_2 v_2 w_2 R_2 u_2$, $\ldots$ , $u_1 P_{2n-6} v_{2n-6} w_{2n-6} R_{2n-6} u_2$, $u_1 u_1^{\prime} u_2$, $u_1 u_2^{\prime} u_2$. Thus, $\kappa(\{u_1,u_2\})\geq 2n-4$.
		
		In conclusion, $\kappa(\widehat{BS_n})=2n-4$. $\square$
	\end{lemma}
	
	
	\begin{figure}[H]
		\centering
		\begin{subfigure}[b]{0.43\textwidth}
			\includegraphics[width=\textwidth]{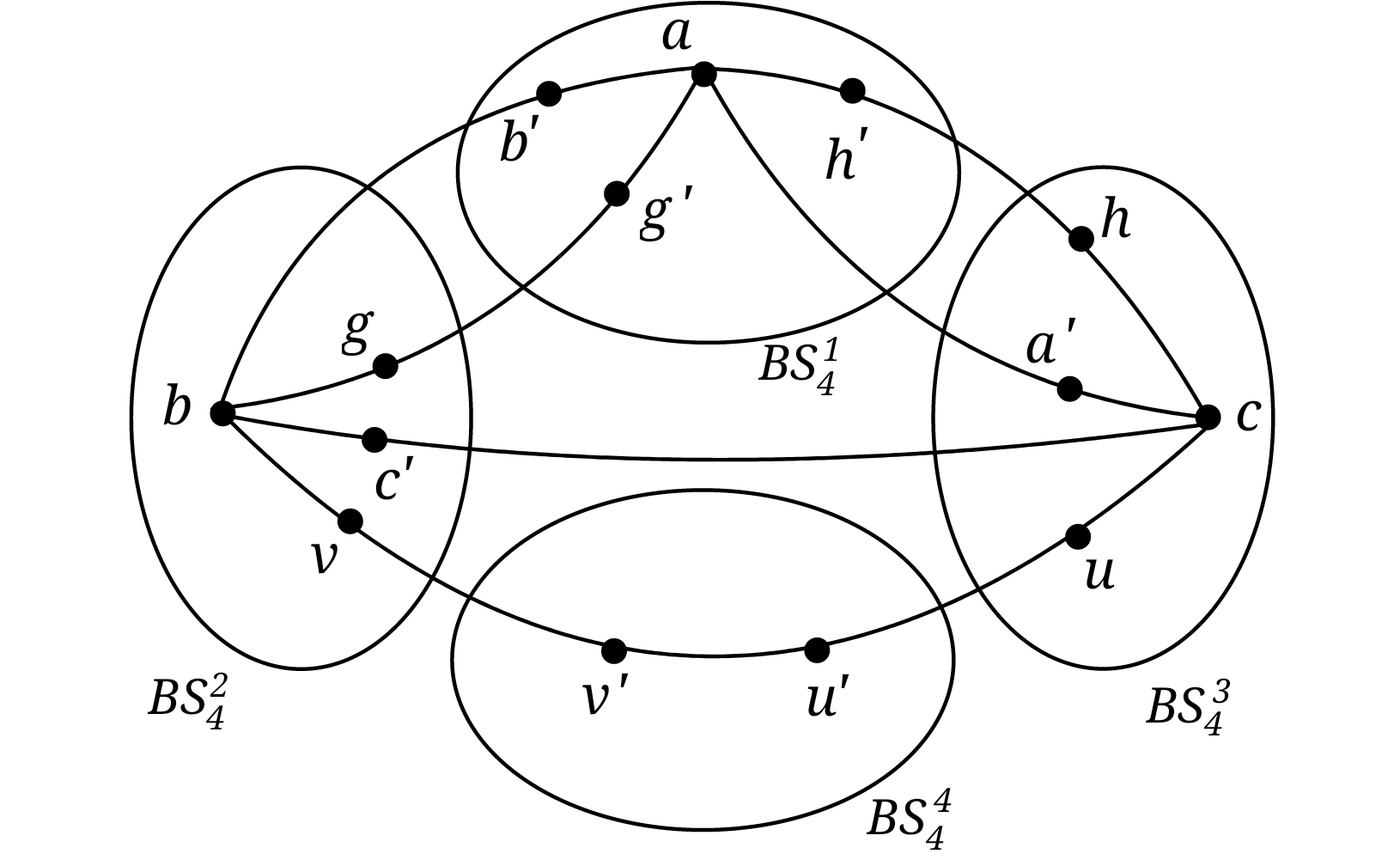} 
			\subcaption{Subcase 3.1}
			
		\end{subfigure}
		\hspace{35pt}
		\begin{subfigure}[b]{0.43\textwidth} 
			\includegraphics[width=\textwidth]{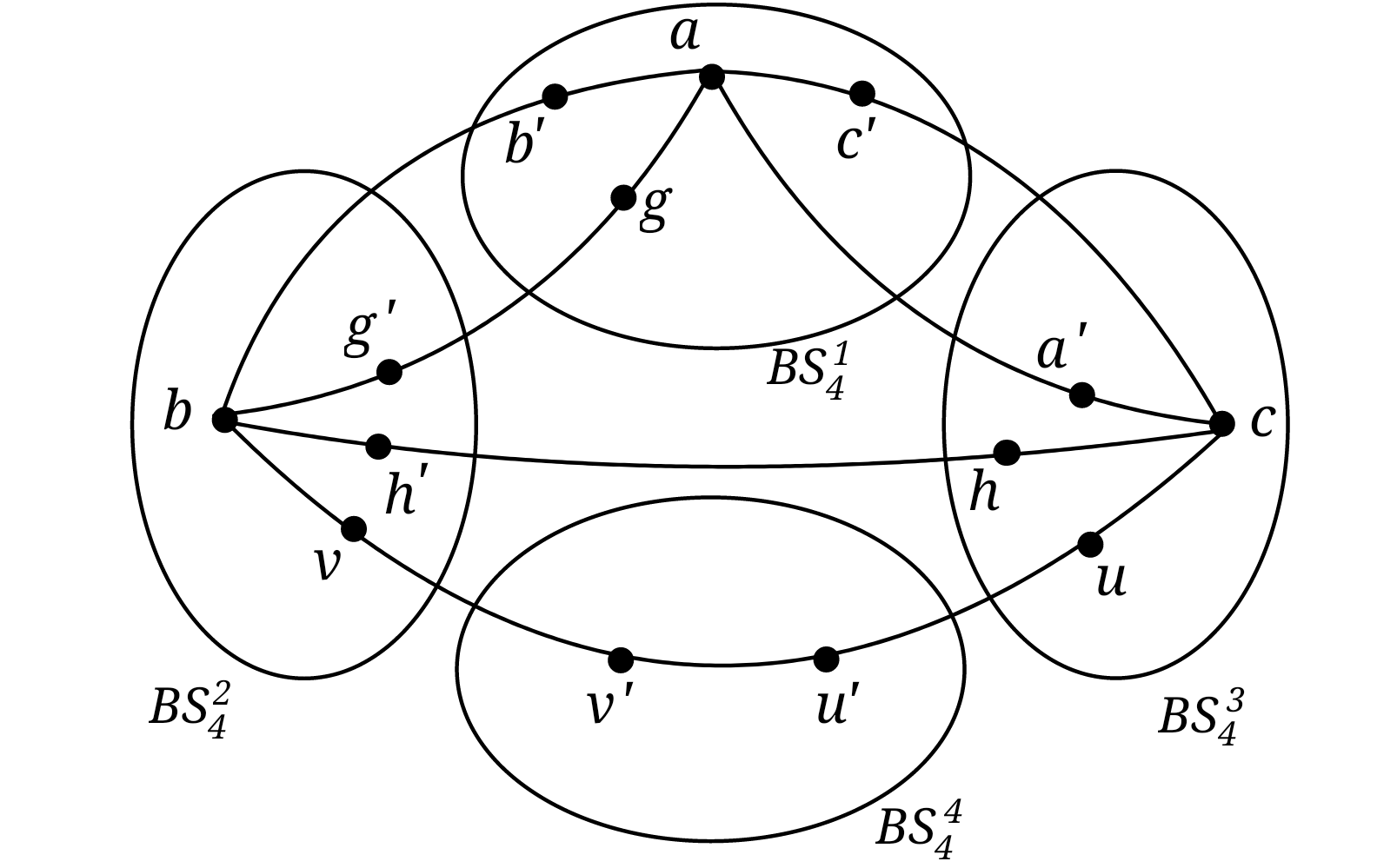} 
			\subcaption{Subcase 3.2.1}
			
		\end{subfigure}
		\caption{The illustrations of Lemma 2.6}
		\label{fig3}
	\end{figure}
	
	\begin{lemma}\label{lem2.6} 
		For any subset $T = \{a, b, c\}$ of $V(BS_4)$, there exist six internally disjoint paths: two $(a, b)$-paths, two $(b, c)$-paths and two $(a, c)$-paths, with none of their internal vertices belonging to $T$.
		\\
		\\
		\noindent\textbf{Proof.}
		We will prove the assertion by considering the following three cases.
		
		\textbf{Case 1:} $a,b$ and $c$ belong to the same copy.
		
		Without loss of generality, suppose that $a,b,c\in BS_4^1\cong BS_3$. By the symmetry of $BS_3$, we only need to consider the following two alternative cases: 
		
		(i) The distance between each pair of the vertices $a$, $b$, and $c$ is $2$. 
		
		(ii) The distances between any two of $a$, $b$ and $c$ are $1$, $1$ and $2$, respectively.
		
		Referring to Fig. \ref{fig1} (b), one can easily check that the assertion holds in the above two cases.
		
		\textbf{Case 2:} $a,b$ and $c$ belong to two different copies.
		
		Without loss of generality, suppose that $a,b\in BS_4^1$ and $c\in BS_4^2$. By Lemma \ref{lem2.1} (i), we have $\kappa(BS_{4}^1)=2n-5=3$, and thus there exist $3$ internally disjoint $(a,b)$-paths in $BS_4^1$. Let $Y=\{a^+,a^-,b^+,b^-\}$. By Lemma \ref{lem2.5}, we have  $\kappa(BS_4\setminus BS_4^1)=4$. Since $|Y|=4$, by Lemma \ref{lem2.3}, there exists a $4$-fan in $BS_4\setminus BS_4^1$ from $c$ to $Y$. Consequently, we can obtain $2$ $(a,c)$-paths: $R_1=aa^+c,~R_2=aa^-c$, and $2$ $(b,c)$-paths: $Q_1=bb^+c,~Q_2=bb^-c$,  that is, the desired structure is obtained. Note that it is possible that $c\in Y$, if so, we can still obtain the desired structure.
		
		\textbf{Case 3:} 
		$a,b$ and $c$ belong to three different copies, respectively.
		
		Without loss of generality, suppose that $a\in BS_4^1, b\in BS_4^2$ and $c\in BS_4^3$. Each of $a,b,c$ has at least one outgoing neighbor in $BS_4^1\oplus BS_4^2\oplus BS_4^3$, and these outgoing neighbors are denoted by   $a^\prime,b^\prime,c^\prime$, respectively. Next, we consider the following two cases:
		
		\textbf{Subcase 3.1:} Each $BS_4^i$ $(i=1,2,3)$ contains exactly one vertex in $\{a^\prime,b^\prime,c^\prime\}$.
		
		Without loss of generality, suppose that $a^\prime\in BS_4^3, b^\prime\in BS_4^1$ and $c^\prime\in BS_4^2$. By Lemma \ref{lem2.1} (ii), there exist vertices $g \in V(BS_4^2) \setminus \{c^\prime\}$ and $h \in V(BS_4^3) \setminus \{a^\prime\}$ such that their outgoing neighbors $g^\prime$ and $h^\prime$ belong to $V(BS_4^1) \setminus \{b^\prime\}$ and $V(BS_4^1) \setminus \{b^\prime, g^\prime\}$, respectively. Additionally, there exist vertices $v \in V(BS_4^2) \setminus \{c^\prime, g\}$ and $u \in V(BS_4^3) \setminus \{a^\prime, h\}$ such that their outgoing neighbors $v^\prime$ and $u^\prime$ belong to $V(BS_4^4)$ (see Fig. \ref{fig3} (a)).

		Since $\kappa(BS_3)=3$, by Lemma \ref{lem2.3}, we can construct a $3$-fan which contains three internally disjoint $(a,\{g^\prime, b^\prime, h^\prime\})$-paths in $BS_4^1$. Similarly, for $i = 2,3$, there exists a fan in $BS_4^i$ that also contains internally disjoint paths analogous to those in the $3$-fan of $BS_4^1$. In addition, there exists a $(v^\prime,  u^\prime)$-path in $BS_4^4$. In summary, the desired structure is obtained.
		
		\textbf{Subcase 3.2:} One of $BS_4^i$ $(i=1,2,3)$ contains exactly two vertices in $\{a^\prime,b^\prime,c^\prime\}$.
		
		Without loss of generality, suppose that $a^\prime\in BS_4^3$ and $ b^\prime,c^\prime\in BS_4^1$. 
		
		\textbf{Subcase 3.2.1:} $c^\prime\neq b^\prime$.
		
		By Lemma \ref{lem2.1} (ii), there exist vertices $g \in V(BS_4^1) \setminus \{b^\prime,c^\prime\}$ and $h \in V(BS_4^3) \setminus \{a^\prime\}$ such that their outgoing neighbors $g^\prime$ and $h^\prime$ belong to $V(BS_4^2)$ and $V(BS_4^2) \setminus \{ g^\prime\}$, respectively. Additionally, there exist vertices $v \in V(BS_4^2) \setminus \{h^\prime, g^\prime\}$ and $u \in V(BS_4^3) \setminus \{a^\prime, h\}$ such that their outgoing neighbors $v^\prime$ and $u^\prime$ belong to $V(BS_4^4)$ (see Fig. \ref{fig3} (b)).
		
		Since $\kappa(BS_4)=3$, by Lemma \ref{lem2.3}, we can construct a $3$-fan from $a$ to $\{b^\prime, c^\prime, g\}$ in $BS_4^1$. Similarly, there exists a fan in $BS_4^i$ $(i=2,3)$. In addition, there exists a $(v^\prime,  u^\prime)$-path in $BS_4^4$.
		In summary, the desired structure is obtained.
		
		\textbf{Subcase 3.2.2:} $c^\prime=b^\prime$.
		
		Without loss of generality, suppose that $c^\prime=b^\prime=c^-=b^+$, in this case we have $c=(2413),~ b=(1432), ~c^-=b^+=(2431)$. One can check the conclusion by simply discussing the cases of $a$.  $\square$
	\end{lemma}

	\begin{lemma}\label{lem2.7} (\cite{Zhu23})
		For any $k$-regular connected graph $G$,  $\pi_3\left(G\right)\leq\lfloor\frac{3k-r}4\rfloor$,
		where $r=max \left \{ \left | CN_{G}( u,v,w) \right| \mid\{u,v,w\}\subseteq V({G})\right\}.$
		
	\end{lemma}
	
	\section{A structure connecting any three vertices in $V(BS_n)$}\label{sec3}
	In this section, considering the parity of $n$, we show that $BS_n$ has a structure connecting any three vertices in $V(BS_n)$.
	
	\begin{theorem}\label{the3.1} 
		Let $T=\{a,b,c\}$ be an arbitrary subset of $V(BS_n)$ with $n\geq4$. Then the following properties hold:
		
		(1) For $n=2k\,(k\geq2)$, the structure contained in $BS_n$ is depicted in Fig. \ref{fig4} (a): there are $6k-6$ internally disjoint paths: $2k-2$ $(a,b)$-paths, $2k-2$ $(b,c)$-paths and $2k-2$ $(a,c)$-paths, with none of their internal vertices belonging to $T$.
		
		(2) For $n=2k+1\,(k\geq2)$, the structure contained in $BS_n$ is depicted in Fig. \ref{fig4} (b): there are $6k-4$ internally disjoint paths: $2k{-}2$ $(a,b)$-paths, $2k{-}2$ $(b,c)$-paths and $2k$ $(a,c)$-paths, with none of their internal vertices belonging to $T$. Additionally, there exist two vertices $b_1,b_2\in N_{BS_n}(b)$ that do not appear on any of these paths.
		\begin{figure}[H]
			\centering
			\begin{subfigure}[b]{0.35\textwidth} 
				\includegraphics[width=\textwidth]{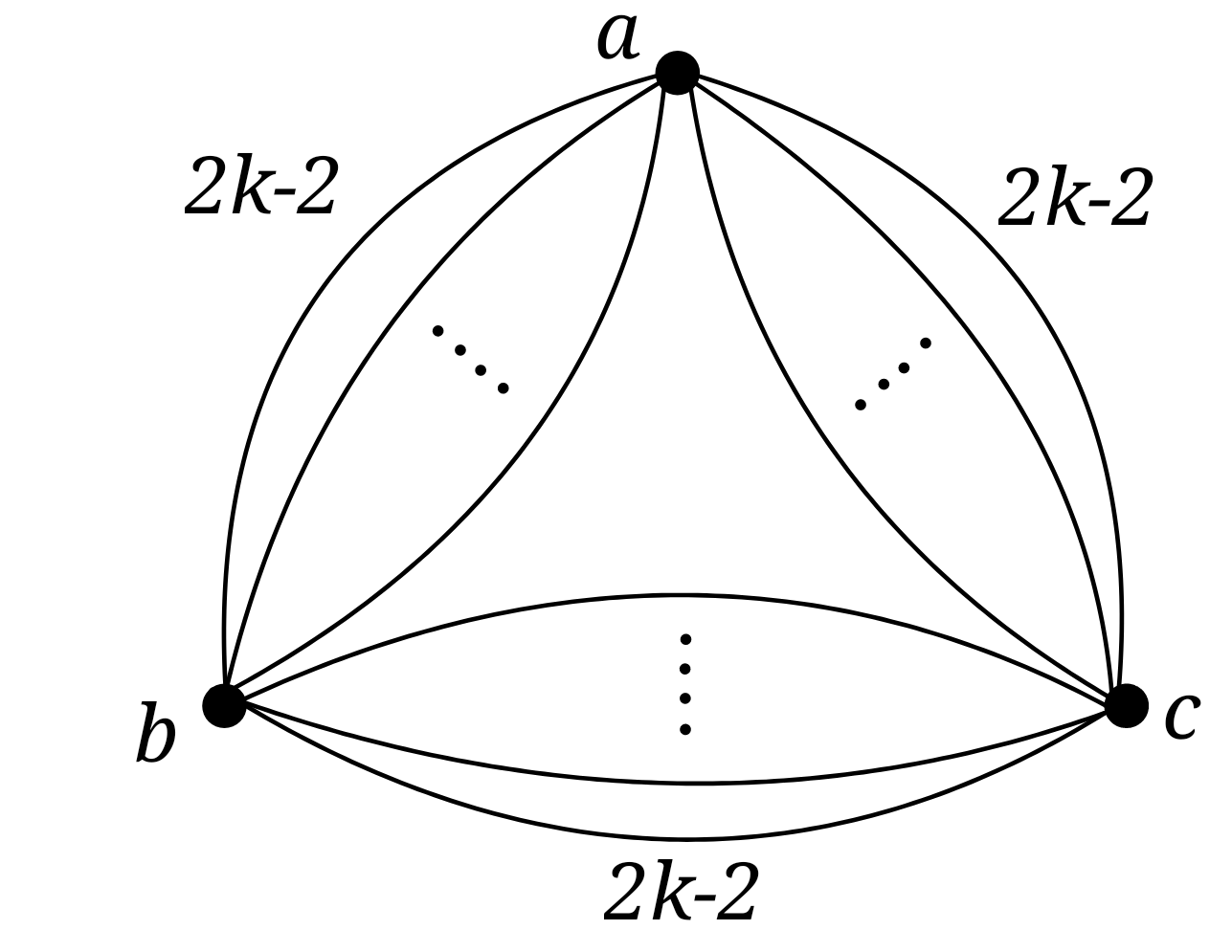} 
				\subcaption{$n=2k\,(k\geq2)$}
				
			\end{subfigure}
			\hspace{35pt}
			\begin{subfigure}[b]{0.35\textwidth} 
				\includegraphics[width=\textwidth]{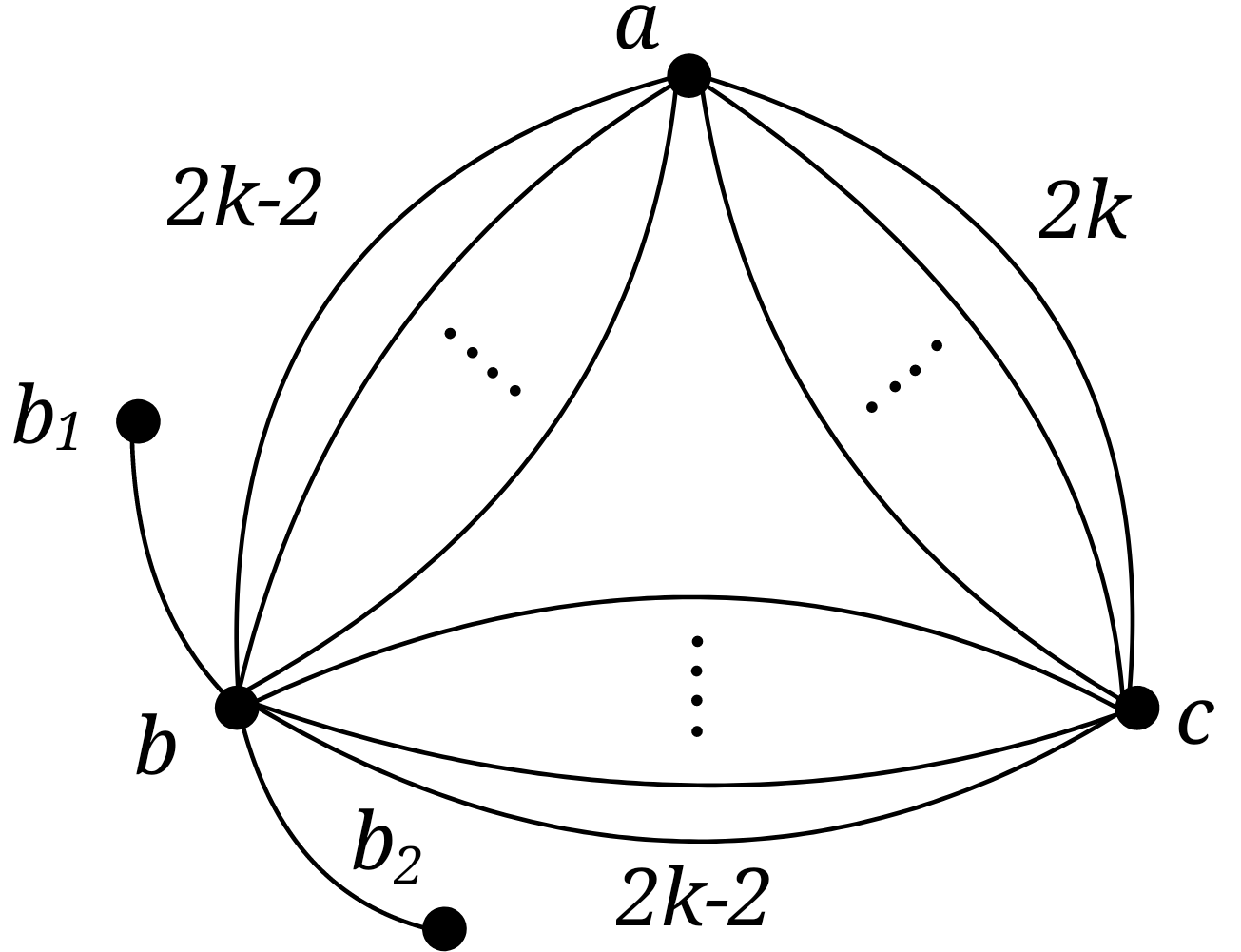} 
				\subcaption{$n=2k+1\,(k\geq2)$}
				
			\end{subfigure}
			\caption{A structure connecting any three vertices in $V(BS_n)$}
			\label{fig4}
		\end{figure}

		\noindent\textbf{Proof.}
		For $n=4$, by Lemma \ref{lem2.6}, the desired structure is obtained.
		Next, we assume that $n\geq 5$, let $BS_n=BS_n^1\oplus BS_n^2\oplus\cdots\oplus BS_n^n$.
		
		We consider the following two parts, which are proved by induction on $n$:
		
		\textbf{Part 1:} Assume that the result holds for $n-1=2k~(k\geq2)$, next we prove that it also holds for $n=2k+1 ~(k\geq2)$.
		
		\textbf{Case 1:} $a,b$ and $c$ belong to the same copy.
		
		Without loss of generality, assume that $a,b,c\in BS_n^1$.  Since $BS_n^1\cong BS_{n-1}$, by induction, the structure depicted in Fig. \ref{fig4} (a) is contained within $BS_n^1$, there are $6k-6$ internally disjoint paths: $2k-2~(a,b)$-paths, $2k-2~(b,c)$-paths and $2k-2~(a,c)$-paths, with none of their internal vertices belonging to $T$.
		
		Let $X=\{a^+,a^-\}$ and $Y=\{c^+,c^-\}$.
		By Lemma \ref{lem2.5}, we have  $\kappa(BS_n\setminus BS_n^1)=2n-4=4k-2$, and thus $\kappa(BS_n\setminus BS_n^1-\{b^+,b^-\})=(2n-4)-2=4k-4>2$. Then, by Lemma \ref{lem2.4}, there are $2$ pairwise disjoint $(X,Y)$-paths in $BS_n\setminus BS_n^1-\{b^+,b^-\}$. Therefore, we can obtain $2~(a,c)$-paths.
		
		In summary, we have obtained $6k-4$ internally disjoint paths: $2k{-}2~(a,b)$-paths, $2k
		{-}2~(b,c)$-paths and $2k~(a,c)$-paths, with none of their internal vertices belonging to $T$. Additionally, let $b_1=b^+,\,b_2=b^-$.
		Then, the desired structure is obtained.
		\begin{figure}[H]
			\centering
			\begin{subfigure}[b]{0.3\textwidth}
				\includegraphics[width=\textwidth]{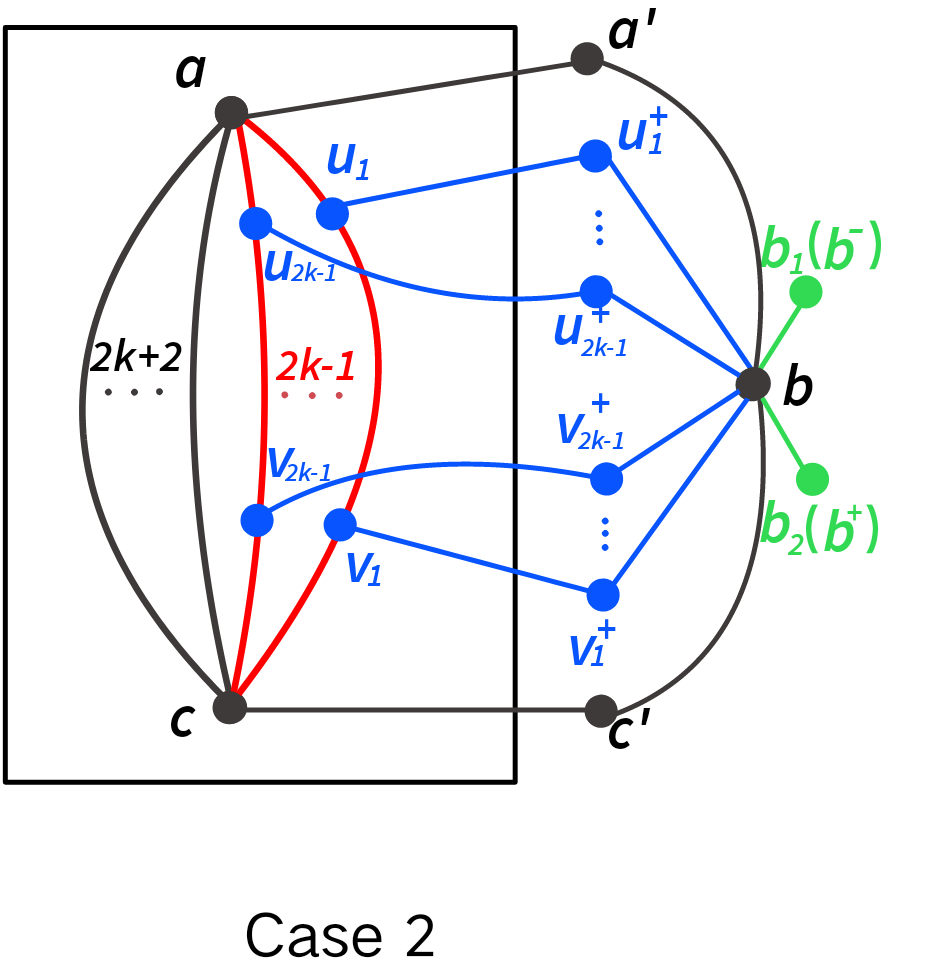}
				\subcaption{Case 2} 
				
			\end{subfigure}
			\hfill
			\begin{subfigure}[b]{0.3\textwidth}
				\includegraphics[width=\textwidth]{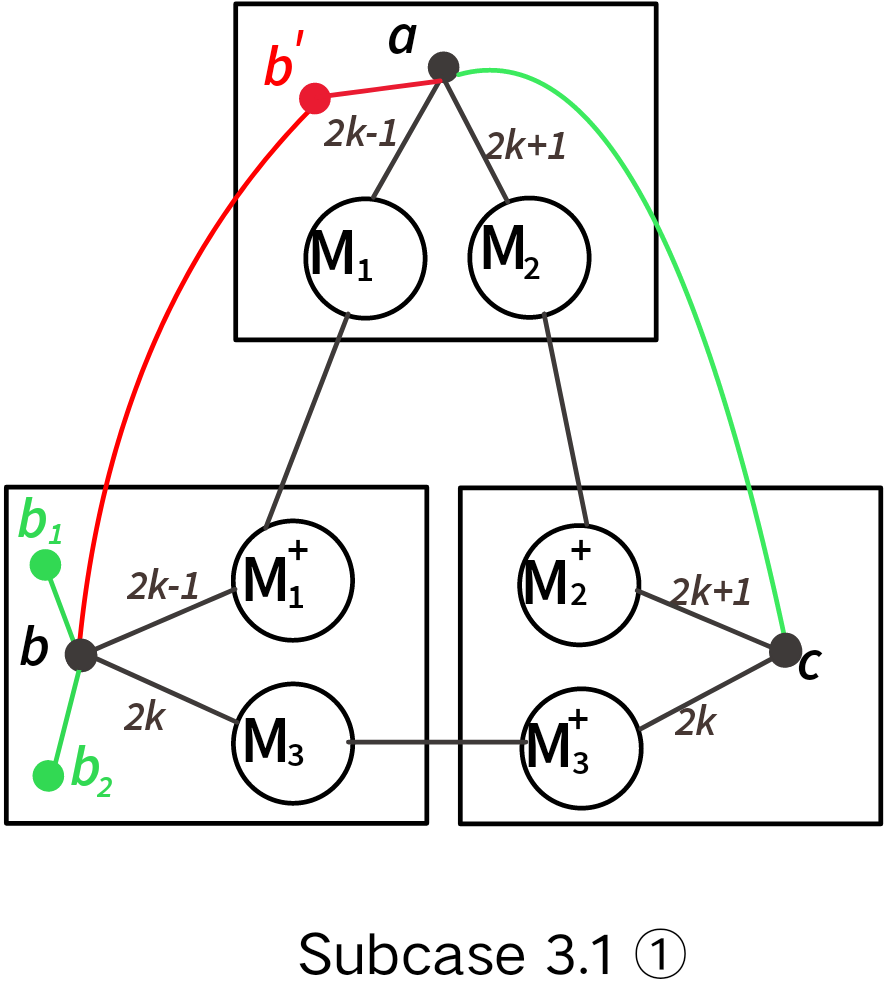}
				\subcaption{Subcase 3.1 \textcircled{1}}
				
			\end{subfigure}
			\hfill
			\begin{subfigure}[b]{0.3\textwidth}
				\includegraphics[width=\textwidth]{2312}
				\subcaption{Subcase 3.1 \textcircled{2}}
				
			\end{subfigure}
			
			\vspace{10pt}
			
			\begin{subfigure}[b]{0.24\textwidth}
				\includegraphics[width=\textwidth]{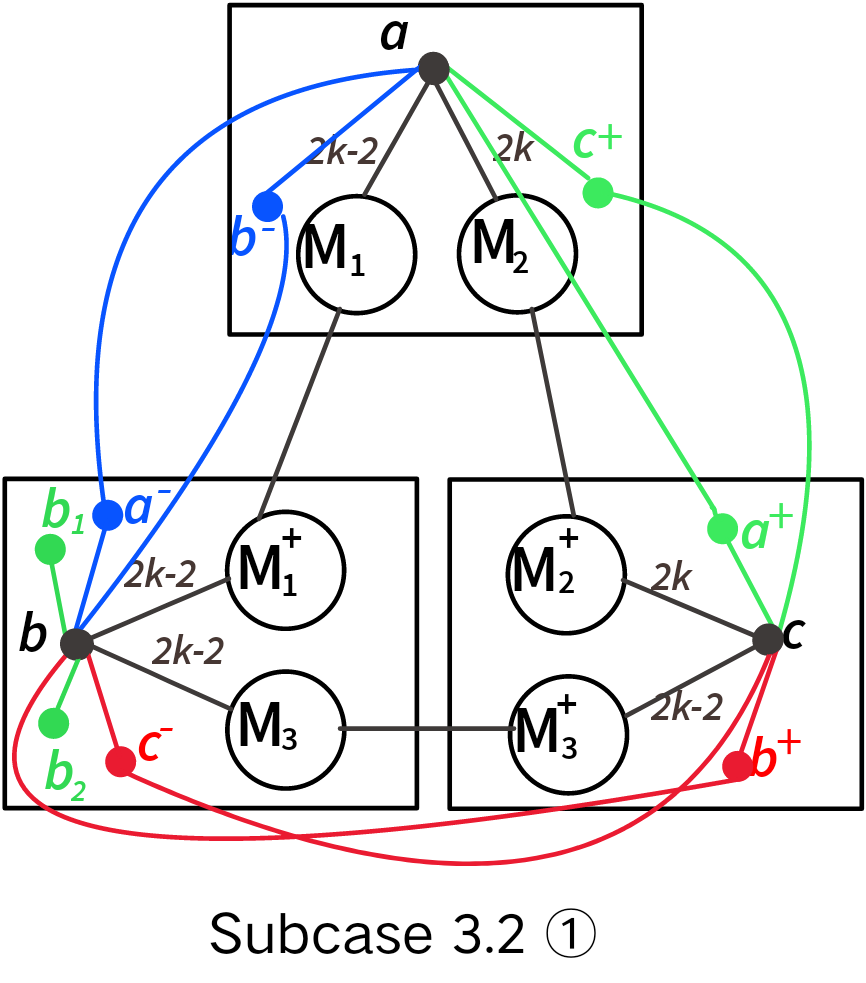}
				\subcaption{Subcase 3.2 \textcircled{1}}
			\end{subfigure}
			\hfill
			\begin{subfigure}[b]{0.24\textwidth}
				\includegraphics[width=\textwidth]{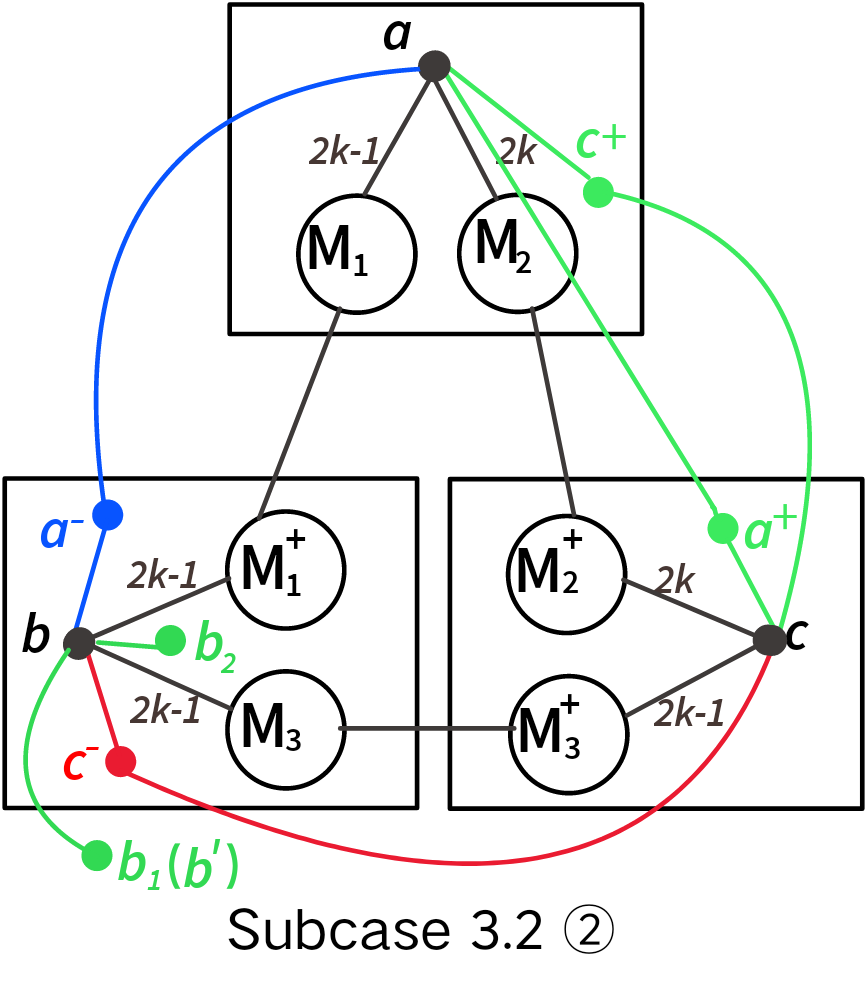}
				\subcaption{Subcase 3.2 \textcircled{2}}
				
			\end{subfigure}
			\hfill
			\begin{subfigure}[b]{0.24\textwidth}
				\includegraphics[width=\textwidth]{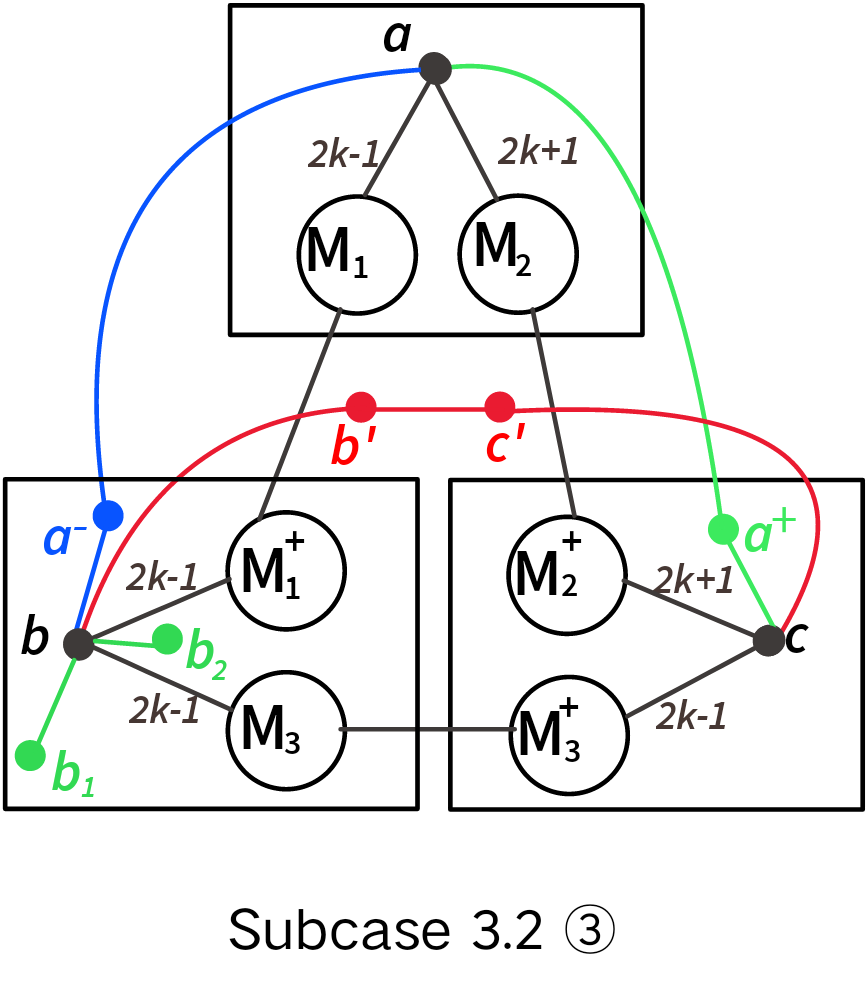}
				\subcaption{Subcase 3.2 \textcircled{3}}
				
			\end{subfigure}
			\hfill
			\begin{subfigure}[b]{0.24\textwidth}
				\includegraphics[width=\textwidth]{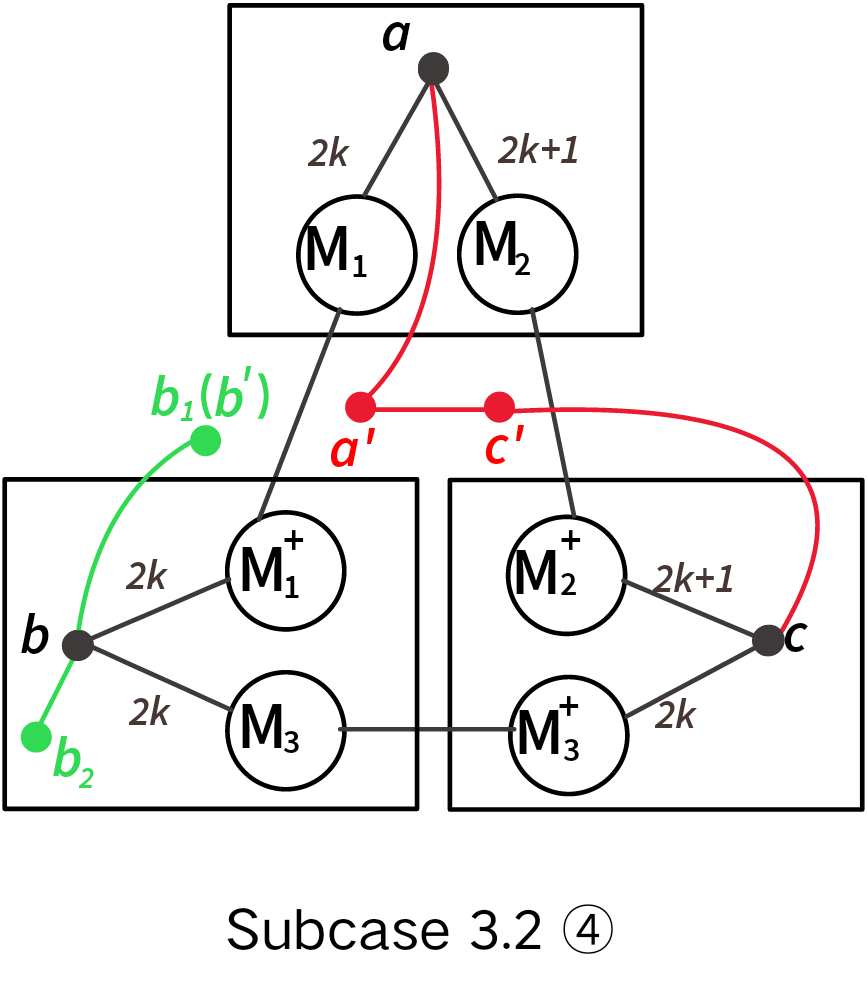}
				\subcaption{Subcase 3.2 \textcircled{4}}
				
			\end{subfigure}
			
			\caption{The illustrations of Part 1}
			\label{fig5}
		\end{figure}
		\textbf{Case 2:}  $a,c$ and $b$ belong to two different copies.
		
		Without loss of generality, assume that $a,c\in BS_n^1$ and $b\in BS_n^2$.
		
		Since $\kappa(BS_n^1)=2n-5=4k-3$, there exist $4k{-}3$ internally disjoint $(a,c)$-paths in $BS_n^1$. By Lemma \ref{lem2.1} (iii), we have  $|CN_{BS_n}(a,c)|\leq3$. Considering that $a,c$ might be adjacent, there are at most $4~(a,c)$-paths of length not exceeding  $2$, and thus there are at least $(4k-3)-4$ $(a,c)$-paths of length at least $3$.
		Since $(2k-3)\leq (4k-3)-4$, we can choose $2k{-}3~(a,c)$-paths $P_1,P_2\ldots,P_{2k-3}$ of length at least $3$. 
		
		Referring to Fig. \ref{fig5} (a), for each $i\in[2k-3]$, 
		there are vertices $u_i,v_i\in V(P_i)$ such that $\{(a,u_i),(c,v_i)\}\subseteq E(P_i)$ and $u_i\neq v_i$. Let $M=\{u_1,u_2,\cdots,u_{2k-3},v_1,v_2,\ldots,v_{2k-3}\}$, $M^+=\{u_1^+,u_2^+,\cdots,u_{2k-3}^+,v_1^+,v_2^+,\ldots,v_{2k-3}^+\}$.
		Respectively, $a$ and $c$ have at least one outgoing neighbor $a^\prime$ and $c^\prime$ such that $\{a^\prime ,c^\prime\}\cap \{b\}=\emptyset$.
		Let $X=M^+\cup\{a^\prime ,c^\prime\}$ and $b_1,b_2\in N_{(BS_n\setminus BS_n^1)-X}(b)$. Clearly, $X\subseteq V((BS_n\setminus BS_n^1)-\{b_1,b_2\})$ and $|X|=4k-4$. Since $\kappa((BS_n\setminus BS_n^1)-\{b_1,b_2\})=2n-4-2=4k-4$, by Lemma \ref{lem2.3}, there exists a family of $4k-4$ internally disjoint $(b,X)$-paths whose terminal vertices are distinct in $X$.
		
		Now we have obtained $2k-2$ $(c,b)$-paths: $R_i=cv_iv_i^+b\,(i\in[2k-3])$, $R_{2k-2}=cc^\prime b$, and $2k-2$ $(a,b)$-paths: $Q_i=au_iu_i^+b\,(i\in[2k-3])$, $Q_{2k-2}=aa^\prime b$. Moreover,  we can get $2k$ $(a,c)$-paths that are distinct from $P_i\,(i\in[2k-3])$. Thus, we can obtain $6k-4$ internally disjoint paths and two vertices $b_1,b_2$ that do not appear on any of these paths.

		\textbf{Case 3:}
		$a,b$ and $c$ belong to three different copies, respectively.
		
		Without loss of generality, assume that $a\in BS_n^1, b\in BS_n^2$ and $c\in BS_n^3$.
		
		Let $H_1=\{v \,|\, v\in V(BS_n^1),v\neq a,v^+\in  V(BS_n^2),\{v^+,v^-\}\cap \{b\}=\emptyset\}$,
		
		~~~~~$H_2=\{v \,|\, v\in V(BS_n^1),v\neq a,v^+\in  V(BS_n^3),\{v^+,v^-\}\cap \{c\}=\emptyset\}$,
		
		~~~~~$H_3=\{v \,|\, v\in V(BS_n^2),v\neq b,v^+\in  V(BS_n^3),\{v^+,v^-\}\cap \{c\}=\emptyset\}$.
		
		Clearly, $H_i \,(i=1,2,3)$ are pairwise disjoint with 	$|H_i|\geq(n-2)!-2\,\textgreater \,2k-1$.

		\textbf{Subcase 3.1:} There exists at least one pair of adjacent vertices in $T$.
		
		Without loss of generality, assume that $a$ and $c$ are adjacent. We consider the following two cases.
		
		\textcircled{1} $b$ has no outgoing neighbor in $BS_n^1\oplus BS_n^3$. 
		
		Referring to Fig. \ref{fig5} (b),  we can find $M_j\subseteq H_j$ such that $|M_j|=2k-2$ $(j=1,3)$ and $M_2\subseteq H_2$ such that $|M_2|=2k-1$ . Let $M_i=\{v_i^1,v_i^2,\ldots,v_i^{|M_i|}\}$ and $M_i^+=\{v^+ | v\in M_i\}$  $(i=1,2,3)$.
		Then, $M_1^+ \subseteq V(BS_n^2)$ and $ M_2^+,M_3^+ \subseteq V(BS_n^3)$.
		
		Since $|M_1\cup M_2|=4k-3$ and  $\kappa(BS_n^1)=2n-5=4k-3$, by Lemma \ref{lem2.3}, there are $4k-3$ internally disjoint $(a,M_1\cup M_2)$-paths $\widetilde{P_1},\widetilde{P_2},\ldots ,\widetilde{P_{2k-2}},$ $\widetilde{R_1},\widetilde{R_2},\ldots, \widetilde{R_{2k-1}}$ such that the terminal vertex of $\widetilde{P_i}$ $(\widetilde{R_j})$ is $v_1^i\,(v_2^j)$ for every $i\in[2k{-}2]$ $(j\in[2k{-}1])$.
		Similarly, there are $4k-4$ internally disjoint $(b,M_1^+ \cup M_3)$-paths $\widehat{P_1},\widehat{P_2},\ldots,\widehat{P_{2k-2}},\widetilde{Q_1}$, $\widetilde{Q_2},\ldots,\widetilde{Q_{2k-2}}$ such that the terminal vertex of $\widehat{P_i}$ $(\widetilde{Q_i})$ is $(v_1^i)^+(v_3^i)$ for every $i\in[2k-2]$. Also, there are  
		$4k-3$ internally disjoint $(c,M_2^+ \cup M_3^+)$-paths $\widehat{R_1},\widehat{R_2},\ldots,\widehat{R_{2k-1}},\widehat{Q_1},\widehat{Q_2},\ldots,\widehat{Q_{2k-2}}$ such that the terminal vertex of $\widehat{R_i}$ $(\widehat{Q_j})$ is $(v_1^i)^+((v_3^j)^+)$ for every $i\in[2k-1]$ $(j\in[2k-2])$.
		
		Let $P_i=a \widetilde{P_i} v_1^i (v_1^i)^+ \widehat{P_i} b$ $(i\in [2k-2])$,  $R_i=a \widetilde{R_i} v_2^i (v_2^i)^+ \widehat{R_i} c$ $(i\in [2k-1])$, $R_{2k}=ac$ and $Q_i=b \widetilde{Q_i} v_3^i (v_3^i)^+ \widehat{Q_i} c$ $(i\in [2k-2])$. We have obtained $6k-4$ internally disjoint paths: $2k-2~(a,b)$-paths $P_i$, $2k-2~(b,c)$-paths $Q_i$ and  $2k~(a,c)$-paths $R_i$, with none of their internal vertices belonging to $T$.
		Additionally, let $b_1=b^+,\,b_2=b^-$. Then, the desired structure is obtained.
		
		\textcircled{2} $b$ has at least one outgoing neighbor in $BS_n^1\oplus BS_n^3$.
		
		Without loss of generality, assume that $b$ has an outgoing neighbor $b^\prime$ in $BS_n^1$.
		Referring to Fig. \ref{fig5} (c), 
		similarly to the proof of Subcase 3.1 \textcircled{1}, we can find $M_1\subseteq H_1$ such that $|M_1|=2k-3$, $M_2\subseteq H_2$ such that $|M_2|=2k-1$ and $M_3\subseteq H_3$ such that $|M_3|=2k-2$. Let $b_1,b_2\in N_{BS_n^2-(M_1^+\cup M_3)}(b)$.
		
		Since $\kappa(BS_{n}^1)=\kappa(BS_{n}^3)=4k-3$ and $|M_1\cup M_2\cup\{b^\prime\}|=|M_2^+\cup M_3^+|=4k-3$, by Lemma \ref{lem2.3}, we can obtain $(4k-3)$ internally disjoint $(a,M_1\cup M_2\cup \{b^\prime\})$-paths in $BS_n^1$ and $(4k-3)$ internally disjoint $(c,M_2^+\cup M_3^+)$-paths in $BS_n^3$. Since $\kappa(BS_{n}^2-\{b_1,b_2\})=4k-5$ and $|M_1^+\cup M_3|=4k-5$, we can obtain $(4k-5)$ internally disjoint $(b,M_1^+\cup M_3)$-paths in $BS_n^2-\{b_1,b_2\}$. By combining these paths, we can construct $6k-4$ internally disjoint paths: $2k-2~(a,b)$-paths, $2k-2~(b,c)$-paths and  $2k~(a,c)$-paths, with none of their internal vertices belonging to $T$. Thus, the desired structure is obtained.

		\textbf{Subcase 3.2:} There are no pairs of adjacent vertices in $T$.
		
		Without loss of generality, we consider the following four cases. 
		
		\textcircled{1} Each vertex in $T$ has no outgoing neighbor in $BS_n\setminus (BS_n^1\oplus BS_n^2\oplus BS_n^3)$.
		
		Referring to Fig. \ref{fig5} (d), without loss of generality, assume that $\{b^-,c^+\}\subseteq V(BS_n^1)$, $ \{a^-,c^-\}\subseteq V(BS_n^2)$, $\{a^+,b^+\}\subseteq V(BS_n^3)$.
		Similarly to the proof of Subcase 3.1 \textcircled{1}, we can find $M_i\subseteq H_i$ such that $|M_i|=2k-4$ $(i=1,3)$ and $M_2\subseteq H_2$ such that $|M_2|=2k-2$. Let $b_1,b_2\in N_{BS_n^2-(M_1^+\cup M_3\cup\{a^-, c^-\})}(b)$.
		
		Since $\kappa(BS_{n}^1)=\kappa(BS_{n}^3)=4k-3$ and $|M_1\cup M_2\cup\{b^-,c^+\}|=|M_2^+\cup M_3^+\cup\{a^+,b^+\}|=4k-4<4k-3$, by Lemma \ref{lem2.3}, we can construct one $(4k-4)$-fan from $a$ to $M_1\cup M_2\cup\{b^-,c^+\}$ in $BS_n^1$ and one $(4k-4)$-fan from $c$ to $M_2^+\cup M_3^+\cup\{a^+,b^+\}$ in $BS_n^3$. Since $\kappa(BS_{n}^2-\{b_1,b_2\})=4k-5$     and $|M_1^+\cup M_3\cup\{a^-, c^-\}|=4k-6<4k-5$, we can construct one $(4k-6)$-fan from $b$ to $M_1^+\cup M_3\cup\{a^-, c^-\}$ in $BS_n^2-\{b_1,b_2\}$. Thus, we can obtain $6k-4$ internally disjoint paths and two vertices $b_1,b_2$ that do not appear on any of these paths. 
		
		\textcircled{2} Exactly one vertex in $T$ has outgoing neighbor in $BS_n\setminus (BS_n^1\oplus BS_n^2\oplus BS_n^3)$. 
		
		Referring to Fig. \ref{fig5} (e), without loss of generality, assume that $b$ has an outgoing neighbor $b^\prime$ in $BS_n\setminus (BS_n^1\oplus BS_n^2\oplus BS_n^3)$, $\{c^+\}\subseteq V(BS_n^1)$, $\{a^-,c^-\}\subseteq V(BS_n^2)$ and $\{a^+\}\subseteq V(BS_n^3)$.  Let $b_1=b^\prime$, $b_2\in N_{BS_n^2-(M_1^+\cup M_3\cup\{a^-,c^-\})}(b)$.   
		Similarly to the proof of Subcase 3.1 \textcircled{1}, we can find $M_i\subseteq H_i$ such that $|M_i|=2k-3~(i=1,3)$ and $M_2\subseteq H_2$ such that $|M_2|=2k-2$.
		
		Since $\kappa(BS_{n}^2-\{b_2\})=4k-4$ and $|M_1^+\cup M_3\cup\{a^-\}|=4k-4$, by Lemma \ref{lem2.3}, we can construct one $(4k-4)$-fan from $b$ to $M_1^+\cup M_3\cup\{a^-,c^-\}$ in $BS_n^2-\{b_2\}$.
		Similarly, we can construct one $(4k-4)$-fan from $a$ to $M_1\cup M_2\cup\{c^+\}$ in $BS_n^1$ and one $(4k-4)$-fan from $c$ to $M_2^+\cup M_3^+\cup\{a^+\}$ in $BS_n^3$. Thus, the desired structure is obtained.
		
		\textcircled{3} Exactly two vertices in $T$ have outgoing neighbor in $BS_n\setminus (BS_n^1\oplus BS_n^2\oplus BS_n^3)$.
		
		Referring to Fig. \ref{fig5} (f), without loss of generality, assume that $\{a^-\}\subseteq V(BS_n^2),~ \{a^+\}\subseteq V(BS_n^3)$ and $b,c$ have outgoing neighbors $b^\prime, c^\prime$ in $BS_n\setminus (BS_n^1\oplus BS_n^2\oplus BS_n^3)$, respectively. Let $b_1,b_2\in N_{BS_n^2-(M_1^+\cup M_3\cup\{a^-\})}(b)$.   
		Similarly to the proof of Subcase 3.1 \textcircled{1},  we can find $M_i\subseteq H_i$ such that $|M_i|=2k-3~(i=1,3)$ and $M_2\subseteq H_2$ such that $|M_2|=2k-1$.
		In addition, there must exist one $(b^\prime,c^\prime)$-path in $BS_n\setminus (BS_n^1\oplus BS_n^2\oplus BS_n^3)$.
		
		Since $\kappa(BS_{n}^2-\{b_1,b_2\})=4k-5$ and $|M_1^+\cup M_3\cup\{a^-\}|=4k-5$, by Lemma \ref{lem2.3}, we can construct one $(4k-5)$-fan from $b$ to $M_1^+\cup M_3\cup\{a^-\}$ in $BS_n^2-\{b_1,b_2\}$.
		Similarly, we can construct one $(4k-4)$-fan from $a$ to $M_1\cup M_2$ in $BS_n^1$ and one $(4k-3)$-fan from $c$ to $M_2^+\cup M_3^+\cup\{a^+\}$ in $BS_n^3$. Thus, the desired structure is obtained.
		
		\textcircled{4} Each vertex in $T$ has at least one outgoing neighbor in $BS_n\setminus (BS_n^1\oplus BS_n^2\oplus BS_n^3)$.
		
		Referring to Fig. \ref{fig5} (g),  without loss of generality, assume that $a,~b$ and $c$ have outgoing neighbors $a^\prime,~ b^\prime$ and $ c^\prime$ in $BS_n\setminus (BS_n^1\oplus BS_n^2\oplus BS_n^3)$, respectively. Since 
		$a^\prime,~ b^\prime$ and $ c^\prime$ cannot all be equal, suppose that $b^\prime\neq a^\prime$ and $b^\prime\neq c^\prime$. Let $b_1=b^\prime$, $b_2\in N_{BS_n^2-(M_1^+\cup M_3)}(b)$.
		Similarly to the proof in Subcase 3.1 \textcircled{1}, we can find $M_i\subseteq H_i$ such that $|M_i|=2k-2~(i=1,3)$ and $M_2\subseteq H_2$ such that $|M_2|=2k-1$. In addition, there must exist one $(a^\prime,c^\prime)$-path in $BS_n\setminus (BS_n^1\oplus BS_n^2\oplus BS_n^3)-\{b^\prime\}$.
		
		Since $\kappa(BS_{n}^2-\{b_2\})=4k-4$ and $|M_1^+\cup M_3|=4k-4$, by Lemma \ref{lem2.3}, we can construct one $(4k-4)$-fan from $b$ to $M_1^+\cup M_3$ in $BS_n^2-\{b_2\}$.
		Similarly, we can construct one $(4k-3)$-fan from $a$ to $M_1\cup M_2$ in $BS_n^1$ and one $(4k-3)$-fan from $c$ to $M_2^+\cup M_3^+$ in $BS_n^3$. Thus, the desired structure is obtained.
		
		\textbf{Part 2:} Assume that the result holds for $n-1=2k+1~(k\geq2)$, next we prove that it also holds for $n=2k+2 ~(k\geq2)$.
		
		\textbf{Case 1:} $a,b$ and $c$ belong to the same copy.
		
		Without loss of generality, assume that $a,b,c\in BS_n^1$. Since $BS_n^1\cong BS_{n-1}$, by induction, the structure depicted in Fig. \ref{fig4} (b) is contained within $BS_n^1$, there are $6k-4$ internally disjoint paths: $2k{-}2~(a,b)$-paths, $2k{-}2~(b,c)$-paths and $2k~(a,c)$-paths, with none of their internal vertices belonging to $T$.

		Let $X=\{a^+,a^-,c^+,c^-\}$ and $Y=\{b^+,b^-,b_1^+,b_2^+\}$. By Lemma \ref{lem2.5} (ii), we have $\kappa(BS_n\setminus BS_n^1)=2n-4=4k>4$. Then, by Lemma \ref{lem2.4}, there are $4$ pairwise disjoint $(X,Y)$-paths in $BS_n\setminus BS_n^1$. Consequently, we can obtain $2~(a,b)$-paths and $2~(b,c)$-paths.
		
		We have obtained $6k$ internally disjoint paths: $2k~(a,b)$-paths, $2k~(b,c)$-paths and $2k~(a,c)$-paths, with none of their internal vertices belonging to $T$. 
		
		\textbf{Case 2:} $a,b$ and $c$ belong to two different copies.
		
		Referring to Fig. \ref{fig6} (a), without loss of generality, assume that $a,b\in BS_n^1$ and $c\in BS_n^2$.
		
		In much the same way as Case 2 in Part 1, since $\kappa(BS_n^1)=2n-5=4k-1$ and $(2k-1)\leq(4k-1)-4$, we can obtain $6k$ internally disjoint paths: $2k~(a,b)$-paths, $2k~(b,c)$-paths and $2k~(a,c)$-paths, with none of their internal vertices belonging to $T$. 
		
		\begin{figure}[H]
			\centering
			\begin{subfigure}[b]{0.3\textwidth}
				\includegraphics[width=\textwidth]{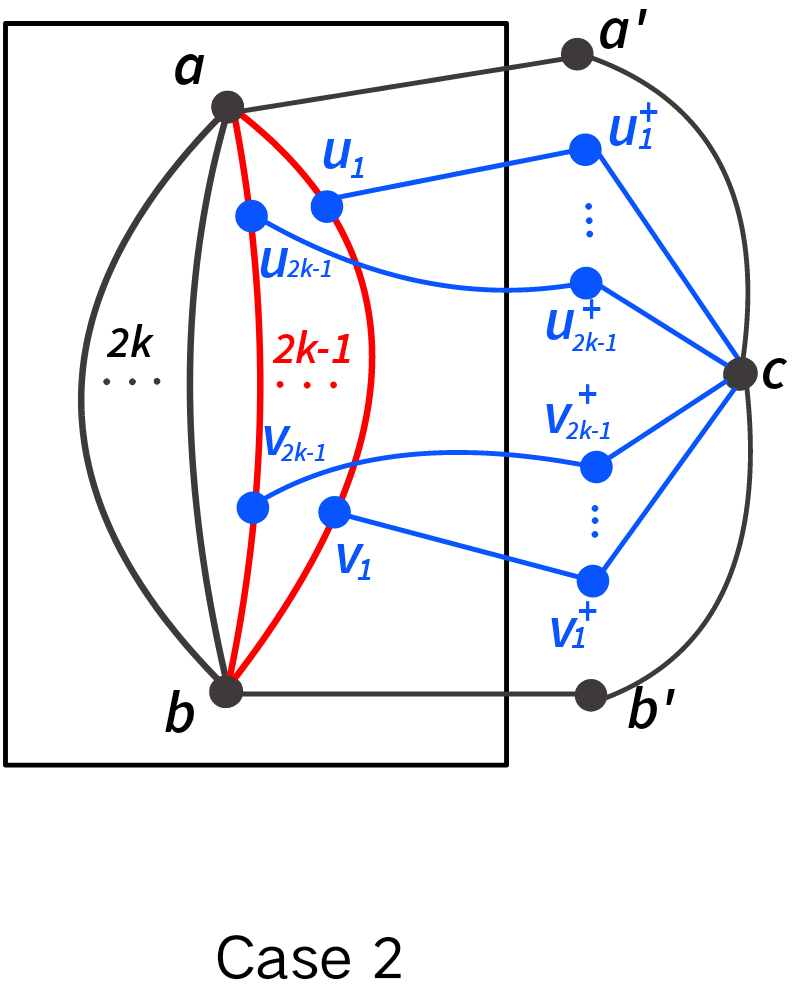}
				\subcaption{Case 2}
			\end{subfigure}
			\hfill
			\begin{subfigure}[b]{0.3\textwidth}
				\includegraphics[width=\textwidth]{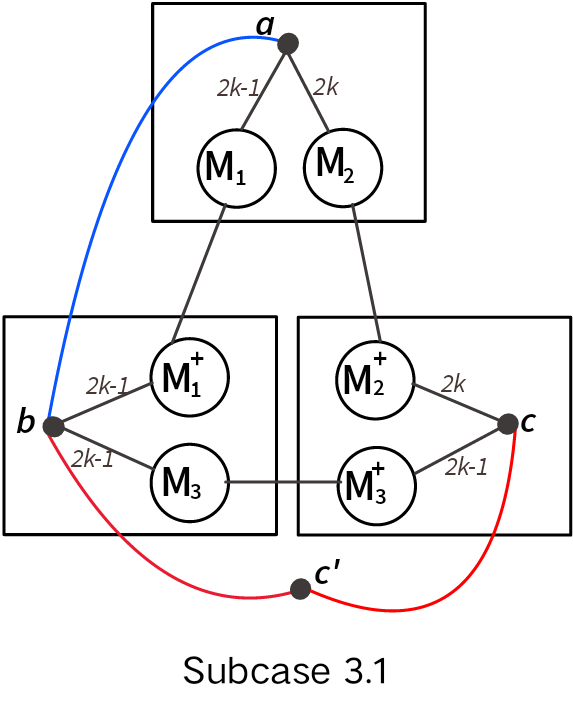}
				\subcaption{Subcase 3.1}
			\end{subfigure}
			\hfill
			\begin{subfigure}[b]{0.3\textwidth}
				\includegraphics[width=\textwidth]{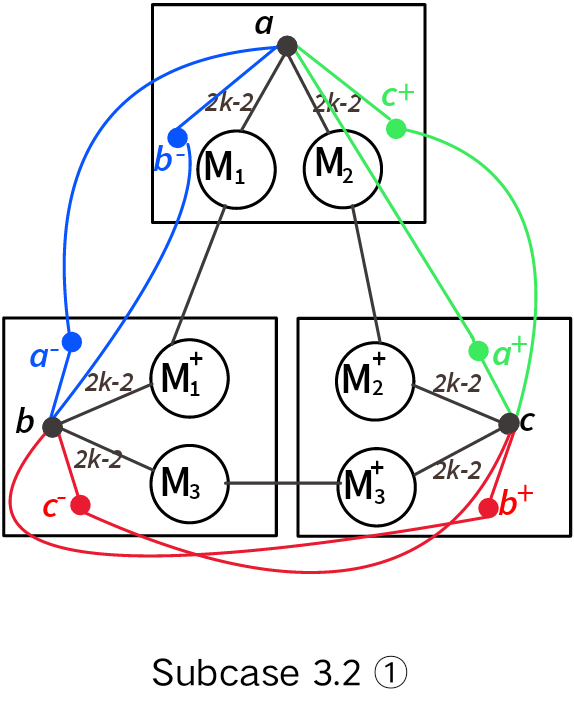}
				\subcaption{Subcase 3.2 \textcircled{1}}
			\end{subfigure}
			
			\vspace{10pt} 
			
			\begin{subfigure}[b]{0.3\textwidth}
				\includegraphics[width=\textwidth]{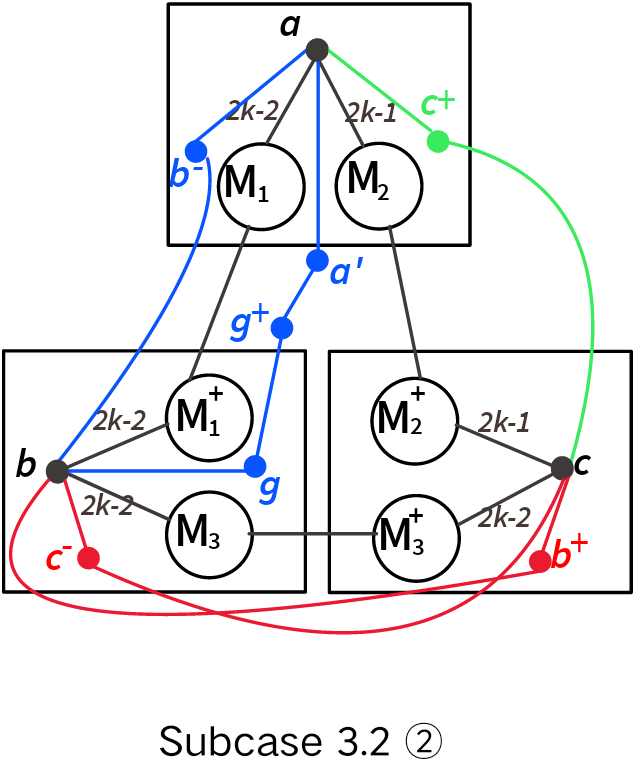}
				\subcaption{Subcase 3.2 \textcircled{2}}
			\end{subfigure}
			\hfill
			\begin{subfigure}[b]{0.3\textwidth}
				\includegraphics[width=\textwidth]{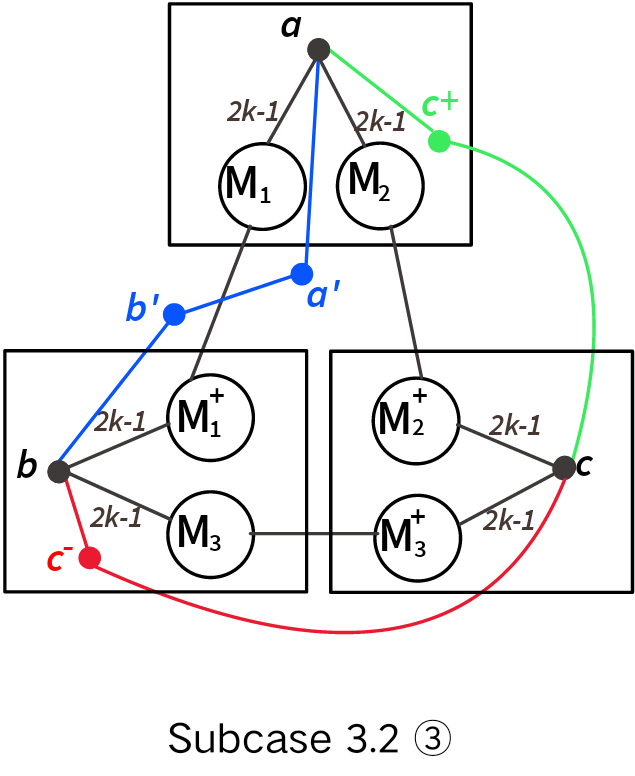}
				\subcaption{Subcase 3.2 \textcircled{3}}
			\end{subfigure}
			\hfill
			\begin{subfigure}[b]{0.3\textwidth}
				\includegraphics[width=\textwidth]{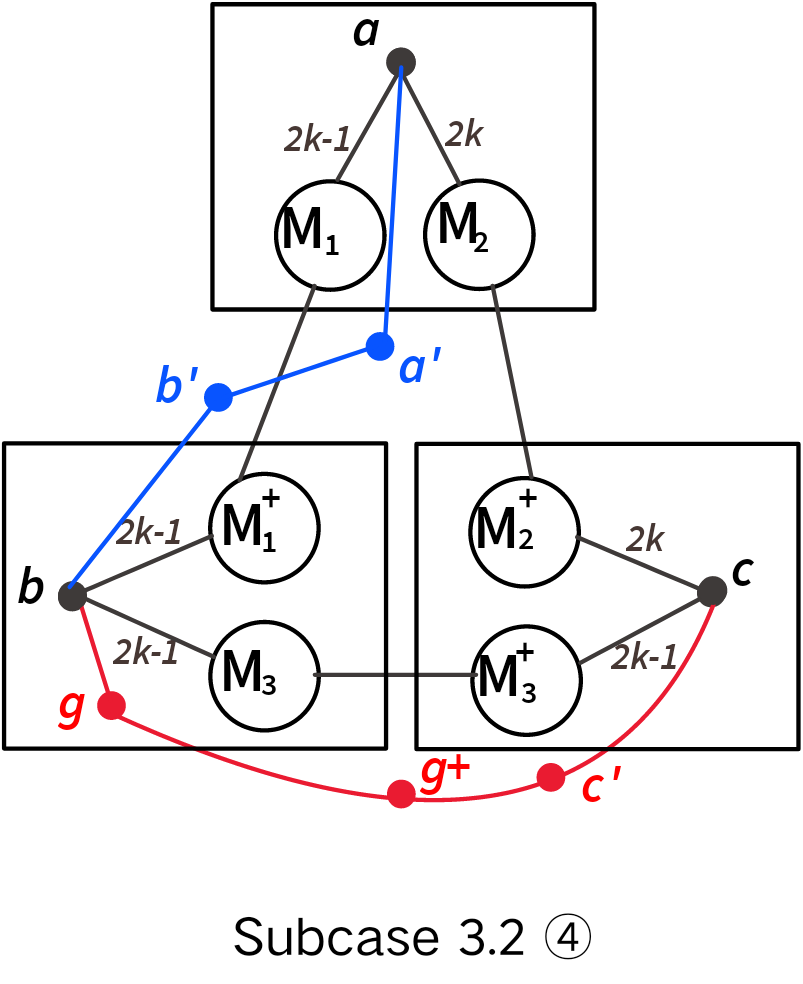}
				\subcaption{Subcase 3.2 \textcircled{4}}
			\end{subfigure}
			
			\caption{The illustrations of Part 2}
			\label{fig6}
		\end{figure}
		\textbf{Case 3:} \label{part2case3}
		$a,b$ and $c$ belong to three different copies, respectively.
		
		Without loss of generality, assume that $a\in BS_n^1, b\in BS_n^2$ and $c\in BS_n^3$.
		
		Let $H_1=\{v \,|\, v\in V(BS_n^1),v\neq a,v^+\in  V(BS_n^2),\{v^+,v^-\}\cap \{b\}=\emptyset\}$,
		
		~~~~~$H_2=\{v \,|\, v\in V(BS_n^1),v\neq a,v^+\in  V(BS_n^3),\{v^+,v^-\}\cap \{c\}=\emptyset\}$,
		
		~~~~~$H_3=\{v \,|\, v\in V(BS_n^2),v\neq b,v^+\in  V(BS_n^3),\{v^+,v^-\}\cap \{c\}=\emptyset\}$.
		
		Clearly, $H_i \,(i=1,2,3)$ are pairwise disjoint with $|H_i|\geq(n-2)!-2>2k$.

		\textbf{Subcase 3.1:} There exists at least one pair of adjacent vertices in $T$.
		
		Referring to Fig. \ref{fig6} (b), without loss of generality, assume that $a$ and $b$ are adjacent, and there must exist an outgoing neighbor $c^\prime$ of $c$ such that $c^\prime \notin V(BS_n^1)$. 
		Similarly to the proof of Subcase 3.1 \textcircled{1} in part 1, we can find $M_2\subseteq H_2$ such that $|M_2|=2k$ and $M_j\subseteq H_j$ such that $|M_j|=2k-1$ $(j=1,3)$. 
		
		Since $|M_1\cup M_2|=|M_2^+\cup M_3^+|=4k-1$ and $\kappa(BS_{n-1})=2n-5=4k-1$, by Lemma \ref{lem2.3}, there are $4k-1$ internally disjoint $(a,M_1\cup M_2)$-paths in $BS_n^1$ and
		$4k-1$ internally disjoint $(c,M_2^+ \cup M_3^+)$-paths in $BS_n^3$. 
		Since $|M_1^+ \cup M_3\cup\{c\prime\}|=4k-1$ and  $\kappa(BS_n\setminus (BS_n^1\oplus BS_n^3))=2n-5=4k-1$, there are $4k-1$ internally disjoint $(b,M_1^+ \cup M_3\cup\{c\prime\})$-paths in $BS_n\setminus (BS_n^1\oplus BS_n^3)$. Then, the desired structure is obtained.
		
		\textbf{Subcase 3.2:} There are no pairs of adjacent vertices in $T$.
		
		Without loss of generality, we consider the following four cases: 
		
		\textcircled{1} Each vertex in $T$ has no outgoing neighbor in $BS_n\setminus (BS_n^1\oplus BS_n^2\oplus BS_n^3)$.
		
		Referring to Fig. \ref{fig6} (c), without loss of generality, assume that  $\{b^-,c^+\}\subseteq V(BS_n^1)$, $\{a^-,c^-\}\subseteq V(BS_n^2)$, $\{a^+,b^+\}\subseteq V(BS_n^3)$. 
		Similarly to the proof of Subcase 3.1 \textcircled{1} in part 1, we can find $M_i\subseteq H_i$ such that $|M_i|=2k-2$ for $i=1,2,3$.  
		
		Since $\kappa(BS_{n-1})=4k-1$ and $|M_1\cup M_2\cup\{b^-, c^+\}|=|M_1^+\cup M_3\cup\{a^-, c^-\}|=|M_2^+\cup M_3^+\cup\{a^+, b^+\}|=4k-2$, by Lemma \ref{lem2.3}, we can construct   one $(4k-2)$-fan from $a$ to $M_1\cup M_2\cup\{b^-, c^+\}$, one $(4k-2)$-fan from  $b$ to $M_1^+\cup M_3\cup\{a^-, c^-\}$  and one $(4k-2)$-fan from $c$ to $M_2^+\cup M_3^+\cup\{a^+, b^+\}$. 
		Then, the desired structure is obtained.

		\textcircled{2} Exactly one vertex in $T$ has outgoing neighbor in $BS_n\setminus (BS_n^1\oplus BS_n^2\oplus BS_n^3)$.
		
		Referring to Fig. \ref{fig6} (d), without loss of generality, assume that $a$ has an outgoing neighbor $a^\prime$ in $BS_n\setminus (BS_n^1\oplus BS_n^2\oplus BS_n^3)$, $\{b^-,c^+\}\subseteq V(BS_n^1),~ \{c^-\}\subseteq V(BS_n^2)$ and $\{b^+\}\subseteq V(BS_n^3)$. 
		Similarly to the proof of Subcase 3.1 \textcircled{1} in part 1, we can find $M_i\subseteq H_i$ such that $|M_i|=2k-2$ $(i=1,3)$ and $M_2\subseteq H_2$ such that $|M_2|=2k-1$. In addition, there must exist a vertex $g\in V(BS_n^2)$ such that $g^+ \in V(BS_n\setminus (BS_n^1\oplus BS_n^2\oplus BS_n^3))$. Since $\kappa(BS_n\setminus (BS_n^1\oplus BS_n^2\oplus BS_n^3)=2n-5\geq1$, there exists one $(g^+,a^\prime)$-path in $BS_n\setminus (BS_n^1\oplus BS_n^2\oplus BS_n^3)$.
		
		Since $\kappa(BS_{n-1})=4k-1$, we can construct one $(4k{-}1)$-fan from $a$ to $M_1\cup M_2\cup\{b^-, c^+\}$, one $(4k{-}2)$-fan from $b$ to $M_1^+\cup M_3\cup\{g, c^-\}$ and one $(4k{-}2)$-fan from $c$ to $M_2^+\cup M_3^+\cup\{b^+\}$. Then, the desired structure is obtained.
		
		\textcircled{3} Exactly two vertices in $T$ have outgoing neighbor in $BS_n\setminus (BS_n^1\oplus BS_n^2\oplus BS_n^3)$.
		
		Referring to Fig. \ref{fig6} (e), without loss of generality, assume that $\{c^+\}\subseteq V(BS_n^1)$, $\{c^-\}\subseteq V(BS_n^2)$ and $a,b$ have outgoing neighbors $a^\prime, b^\prime$ in $BS_n\setminus (BS_n^1\oplus BS_n^2\oplus BS_n^3)$, respectively.  
		Similarly to the proof of Subcase 3.1 \textcircled{1} in part 1, we can find $M_i\subseteq H_i$ such that $|M_i|=2k-1$ for $i=1,2,3$. In addition, there must exist one $(a^\prime,b^\prime)$-path in $BS_n\setminus (BS_n^1\oplus BS_n^2\oplus BS_n^3)$.
		
		Since $\kappa(BS_{n-1})=4k-1$, we can construct one $(4k{-}1)$-fan from $a$ to $M_1\cup M_2\cup\{c^+\}$, one $(4k{-}1)$-fan from $b$ to $M_1^+\cup M_3\cup\{c^-\}$ and one $(4k{-}2)$-fan from $c$ to $M_2^+\cup M_3^+$. Then, the desired structure is obtained.
		
		\textcircled{4} Each vertex in $T$ has at least one outgoing neighbor in $BS_n\setminus (BS_n^1\oplus BS_n^2\oplus BS_n^3)$.
		
		Referring to Fig. \ref{fig6} (f), without loss of generality, assume that $a,~b$ and $c$ have  outgoing neighbors $a^\prime,~ b^\prime$ and $ c^\prime$  in $BS_n\setminus (BS_n^1\oplus BS_n^2\oplus BS_n^3)$, respectively. 
		Similarly to the proof of Subcase 3.1 \textcircled{1} in part 1, we can find $M_i\subseteq H_i$ such that $|M_i|=2k-1$ $(i=1,3)$ and $M_2\subseteq H_2$ such that $|M_2|=2k$. In addition, there must exist $g\in V(BS_n^2)$ such that $g^+ \in V(BS_n\setminus (BS_n^1\oplus BS_n^2\oplus BS_n^3))$. Let $X=\{g^+,b^\prime\}, ~Y=\{a^\prime, c^\prime\}$, there are 2 disjoint $(X,Y)$-paths in $BS_n\setminus$ $(BS_n^1\oplus BS_n^2\oplus BS_n^3)$.
		
		Since $\kappa(BS_{n-1})=4k-1$, we can construct one $(4k{-}1)$-fan from $a$ to $M_1\cup M_2$, one $(4k{-}1)$-fan from $b$ to $M_1^+\cup M_3\cup\{g\}$ and one $(4k{-}1)$-fan from $c$ to $M_2^+\cup M_3^+$. Then, the desired structure is obtained.
		
		Combining the results of Part 1 and Part 2, the proof is completed. $\square$
		
	\end{theorem}
	\section{Main result}\label{sec4}
	
	\begin{theorem}\label{the4.1}
		
		Let $BS_n$ be the n-dimensional bubble-sort star graph with $n\geq3$, then $\pi_3(BS_n)=\lfloor\frac{3n}2\rfloor-3$.
		\\
		\\
		\noindent\textbf{Proof.} 
		For $n=3$, by Lemma \ref{lem2.2} and Lemma \ref{lem2.7}, we have $\pi_3(BS_3)\leq\lfloor\frac{3\times3-3}4\rfloor=1$. On the other hand, by the definition of $3$-path-connectivity, it is easy to see that $\pi_3(BS_3)\geq1$. Thus, $\pi_3(BS_3)=1=\lfloor\frac{3\times3}2\rfloor-3$.
		
		For $n\geq4$, since $BS_n$ is a $(2n-3)$-regular connected graph, by Lemma \ref{lem2.2} and Lemma \ref{lem2.7}, $\pi_3(BS_n)\leq\lfloor\frac{3n}2\rfloor-3$.
		Let $T=\{a,b,c\}$ be an arbitrary subset of $V(BS_n)$. Next, we only need to prove that $\pi_3(BS_n)\geq\lfloor\frac{3n}2\rfloor-3$ by considering the parity of $n$.
		
		\textbf{Case 1:} $n=2k$.
		
		In this case, $\lfloor\frac{3n}2\rfloor-3=3k-3$.
		
		By Theorem \ref{the3.1}, $BS_n$ has the structure that is depicted in Fig. \ref{fig4} (a). There are $2k-2$ $(a,b)$-paths, $2k-2$ $(b,c)$-paths and $2k-2$  $(a,c)$-paths. Moreover, these $6k-6$ paths are internally disjoint, and none of their internal vertices belongs to the set $T$.
		
		We first match $k-1$ paths out of $2k-2~ (a,b)$-paths with $k-1$ paths out of $2k-2~ (a,c)$-paths. 
		Next, we match the remaining $k-1$ paths out of $2k-2~(a,b)$-paths with $k-1$ paths out of $2k-2~(b,c)$-paths. 
		Finally, we match the remaining $k-1$ paths out of $2k-2~ (a,c)$-paths with the remaining $k-1$ paths out of $2k-2~(b,c)$-paths.	
		In this manner, we can find $3k-3$ internally disjoint $T$-paths, and thus $\pi_3(BS_n)\geq\lfloor\frac{3n}2\rfloor-3$.
		
		\textbf{Case 2:} $n=2k+1$.
		
		In this case, $\lfloor\frac{3n}2\rfloor-3=3k-2$.
		
		By Theorem \ref{the3.1}, $BS_n$ has the structure that is depicted in Fig. \ref{fig4} (b). There are $2k{-}2$ $(a,b)$-paths, $2k{-}2$ $(b,c)$-paths and $2k$ $(a,c)$-paths. Moreover, these $6k-4$ paths are internally disjoint, and none of their internal vertices belongs to the set $T$.
		
		We first match $k$ paths out of $2k{-}2~ (a,b)$-paths with $k$ paths out of $2k~ (a,c)$-paths. 
		Next, we match the remaining $k-2$ paths out of $2k{-}2~(a,b)$-paths with $k-2$ paths out of $2k{-}2~(b,c)$-paths. 
		Finally, we match the remaining $k$ paths out of $2k~ (a,c)$-paths with the remaining $k$ paths out of $2k{-}2~(b,c)$-paths. 
		In this manner, we can find $3k-2$ internally disjoint $T$-paths, and thus $\pi_3(BS_n)\geq\lfloor\frac{3n}2\rfloor-3$.

		Due to the arbitrariness of $T=\{a,b,c\}$, we have proved that $\pi_3(BS_n)\geq\lfloor\frac{3n}2\rfloor-3$.
		
		Hence, $\pi_3(BS_n)=\lfloor\frac{3n}2\rfloor-3$. $\square$
		
	\end{theorem}

	
\end{document}